\documentclass{amsart}
\usepackage{amsmath}
\usepackage[mathscr]{euscript}
  \usepackage{paralist}
  \usepackage{graphics} %% add this and next lines if pictures should be in esp format
  \usepackage{epsfig} %For pictures: screened artwork should be set up with an 85 or 100 line screen
 \usepackage[colorlinks=true]{hyperref}
\hypersetup{urlcolor=blue, citecolor=red}

\usepackage{JTstuff}

\newcommand{\dee}{{\rm d}}
\newcommand{\rmi}{\mathrm{i}}

\newcommand{\I}{\mathcal{I}}
\newcommand{\MMM}{\mathscr{M}}

\newcommand{\imagepath}{}
\title[D-bar for EIT with conformal maps]{Enhancing D-bar reconstructions for electrical impedance tomography with conformal maps}

\author{Nuutti Hyv\od nen}
\address{Aalto University, Department of Mathematics and Systems Analysis, P.O. Box 11100, FI-00076 Aalto, Finland}
\email{nuutti.hyvonen@aalto.fi}

\author{Lassi P\"aiv\"arinta} 
\address{Tallinn University of Technology, Department of Mathematics, Ehitajate tee 5, 19086 Tallinn, Estonia} 
\email{lassi.paivarinta@ttu.ee}
\thanks{The work of NH and JT was supported by the Academy of Finland (decision 267789). The work of LP was supported by Estonian government grant PUT1093.}

\author{Janne P.~Tamminen}
\address{Aalto University, Department of Mathematics and Systems Analysis, P.O. Box 11100, FI-00076 Aalto, Finland}
\email{janne.tamminen@aalto.fi}


\begin{document}

\begin{abstract}
We present a few ways of using conformal maps in the reconstruction of two-dimensional conductivities in electrical impedance tomography. First, by utilizing the Riemann mapping theorem, we can transform any simply connected domain of interest to the unit disk where the D-bar method can be implemented most efficiently. In particular, this applies to the open upper half-plane. Second, in the unit disk we may choose a region of interest that is magnified using a suitable M{\"o}bius transform. 
To facilitate the efficient use of conformal maps, we introduce input current patterns that are named \textit{conformally transformed truncated Fourier basis}; in practice, their use corresponds to positioning the available electrodes close to the region of interest. These ideas are numerically tested using simulated continuum data in bounded domains and simulated point electrode data in the half-plane. The connections to practical electrode measurements are also discussed.
\end{abstract}

\maketitle

\bigskip

%\tableofcontents

\section{Introduction}
{\em Electrical impedance tomography} (EIT) is an intensively studied research topic because it is a fruitful mathematical problem and has also applications in medical imaging and non-destructive testing. See \cite{Borcea2002,Cheney1999,Uhlmann2009} for overviews of EIT. To summarize, the goal of EIT is to form an image of the conductivity distribution inside a physical body by using electric measurements on its boundary. This corresponds to a nonlinear, ill-posed inverse problem, which needs to be tackled by some regularized algorithm or by resorting to Bayesian inference. In this work, we use the \textit{D-bar method}, which is an example of a direct reconstruction algorithm that assumes the knowledge of the idealized boundary measurements,~i.e.,~the \textit{Dirichlet-to-Neumann} (D-to-N) map. To be more precise, the D-to-N map is linked to the unknown conductivity by equations found by investigating the inverse problem analytically. 

The D-to-N map is a theoretical object representing all possible potential-to-current measurements on the boundary of the domain of interest. Alternatively, one may consider the \textit{Neumann-to-Dirichlet} (N-to-D) map, broadly speaking the inverse of D-to-N map, which sends the applied boundary current to the measured boundary potential. The N-to-D map is arguably the better way to model real measurements: It is smoothening and thus easier to handle numerically, and the approximative connection between real-world electrode measurements and the N-to-D map has also been analyzed~\cite{Hyvonen2009,Hanke2011}. In particular, the boundary condition between the employed electrodes is in practice {\em always} of the homogeneous Neumann type, which obviously favors the use of the N-to-D map (cf.~\cite{Hyvonen2009}). Due to these reasons, we simulate all boundary measurements employed in this paper by computing (suitably truncated versions of) the N-to-D maps for the considered target configurations, and subsequently form the needed (relative) D-to-N maps only when applying the D-bar method to reconstructing conductivities in the {\em unit disk}. Such an approach is enabled in two spatial dimensions by the use of conformal mappings as explained in what follows.

In this paper, we consider EIT in two-dimensional simply connected domains; observe that such considerations do have practical relevance since a two-dimensional model can be utilized if the imaged three-dimensional object is cylindrically symmetric with homogeneous Neumann conditions at the top and bottom boundaries. The original imaging problem is transformed to the unit disk with the help of a suitable conformal map --- by doing so, we enable the use of our direct reconstruction method in the simplest possible setting. 
Furthermore, in the unit disk we may use M\od bius transforms to `magnify' a {\em region of interest} (ROI). Both of these steps are implemented by modifying the applied boundary current densities so that they account for the stretching introduced by the conformal maps and correspond to the standard (truncated) Fourier boundary current basis for the unit disk where the D-bar algorithm is implemented. We call the transformed domain and conductivity the \textit{virtual domain} and \textit{virtual conductivity}, respectively, as opposed to the {\em true domain} and {\em true conductivity} that correspond to the actual measurements. As demonstrated by our numerical experiments, forming the reconstruction in the virtual domain results in an improved outcome inside the ROI.

Magnifying the ROI by a M\od bius transform is motivated by practical measurement set-ups where the number of available electrodes is restricted, but their positions can be freely chosen: If the electrodes are placed equiangledly on the boundary of the unit disk, their number, loosely speaking, restricts the number of Fourier modes that can be used as input current densities (see~\cite{Hyvonen2009,Hanke2011} and Appendix~\ref{sec:point}). Similarly, if one wants to get more accurate information about a certain ROI inside the imaged object, the pattern of electrodes should arguably be the densest on the boundary section closest to the ROI, enabling the use of current patterns with higher spatial frequencies there. The magnification by a M\od bius transform can thus be thought of as mapping a nonuniform electrode pattern designed for accurate imaging of the ROI onto an equiangled pattern that facilitates a straightforward reconstruction by the D-bar method in the unit disk. We give further evidence for this interpretation by presenting theoretical results on the possibility of approximating a given (continuum) current density by electrode measurements modeled by the {\em point electrode model}~(PEM)~\cite{Hanke2011}. This would also straightforwardly lead to similar results for the {\em complete electrode model}~(CEM)~\cite{Cheng1989} with small electrodes.

On a general level, the idea of using conformal mappings in EIT-related considerations is definitely not new: In \cite{Seagar1984} the M{\"o}bius transforms between disks and the half-plane were already presented in the discussion of \textit{visibilities} and \textit{sensitivities} of computed tomography. This discussion continued as a question of \textit{distinguishability} for EIT in \cite{Isaacson1986,Cheney1992}, using again conformal mappings in \cite{Winkler2014,Garde2016}. Following the ideas presented in \cite{Winkler2014}, it was recently demonstrated in \cite{Hyvonen2016} that the CEM is approximately conformally invariant (with constant contact resistances). M{\"o}bius transforms have also been employed in the implementation of the so-called {\em convex source support} algorithm for detecting anomalies from a single measurement pair of EIT both in a bounded domain~\cite{Hanke2008b} and in a half-plane~\cite{Harhanen2010}. This paper uses the  same transformation properties of conformal maps as the aforementioned articles, but we also introduce and numerically test the concept of virtual domain with the help of the {\em conformally transformed truncated Fourier basis} introduced in \cite{Garde2016}. In particular, to the authors' best knowledge, this work presents the first D-bar reconstructions in a half-space geometry.

Implementing a numerical method for computing conformal maps between planar domains is outside the scope of this manuscript. Therefore, we employ the {\em Schwarz--Christoffel Toolbox}~\cite{Driscoll1996} in MATLAB to compute numerical conformal maps and their derivatives between the unit disk and arbitrary polygons~\cite{Trefethen2002}. On the other hand, M\od bius transforms and their derivatives are trivial to implement numerically.

This text is organized as follows. In Section~\ref{sec:setting}, we introduce our mathematical framework. The main ideas of the D-bar method are recalled in Section~\ref{sec:Dbar}, with special emphasis on its numerical implementation (in the unit disk). Section~\ref{sec:conformal} briefly considers how a transformation induced by a conformal mapping affects the Neumann boundary value problem defining the N-to-D map. 
The numerical results are presented in Section~\ref{sec:numerics}. The treatment of half-space geometry as well as the theoretical results related to the approximation of (relative) N-to-D maps with the help of pointlike electrodes are included in two appendices.

\section{The setting}
\label{sec:setting}
The inverse problem of EIT is mathematically described by the \textit{Calder{\'o}n problem}~\cite{Calder'on1980}. Let $\Om \subset \R^2$ be a simply connected, bounded domain with a Lipschitz boundary and let
$$
\sigma\in L^\infty_+(\Om) \, := \, \{\kappa \in L^\infty(\Om) \ | \ {\rm ess} \inf \kappa > 0\} 
$$ 
be a strictly positive, real-valued, isotropic conductivity. A boundary voltage distribution $g\in H^{1/2}(\DOm)$ creates an electrostatic interior potential $u \in H^1(\Omega)$ that solves the Dirichlet problem
\begin{eqnarray}
\label{eq:dirichlet}
  \left\{ \begin{array}{ll}
  \nabla \cdot (\sigma \nabla u) =0  \quad & \rm{in} \  \Om \label{ody1},  \\[1mm]
  u=g  \quad &\rm{on} \ \partial \Om.
\end{array} \right.
\end{eqnarray}
The resulting current density through the boundary is
\begin{equation}\label{DNdef}
  \Lambda_\sigma g := \sigma\frac{\partial u}{\partial \nu}\big|_{\partial\Om} \in H^{-1/2}_\diamond(\partial \Omega),
\end{equation}
where  $\nu \in L^\infty(\partial \Omega, \R^2)$ is the outward unit normal of $\partial \Omega$, the linear operator $\Lambda_\sigma$ is the D-to-N map and
$$
H_{\diamond}^s(\partial\Om) \, := \, \{ v \in  H_{\diamond}^s(\partial\Om) \ | \ \langle f, 1 \rangle_{\partial \Omega} = 0\}, \qquad s \in \R,
$$ 
is the subspace of $H^s(\partial\Om)$ consisting of those functions/distributions that have vanishing mean. Here and in what follows, $\langle \cdot, \cdot \rangle_{\partial \Omega}: H^s(\partial\Om) \times H^{-s}(\partial\Om) \to \C$ denotes the (bilinear) dual bracket on $\partial \Omega$. Since the boundary potentials $g$ and $g + c$, $c \in \C$, for \eqref{eq:dirichlet} correspond to the interior potentials $u$ and $u + c$, respectively, we may define the D-to-N map as
$$
\Lambda_\sigma: H^{1/2}(\partial\Om)/\C \to H_{\diamond}^{-1/2}(\partial\Om),
$$
where the use of a quotient space reflects the irrelevance of the ground level of potential in modeling the underlying physical phenomenon. It follows from the standard theory of elliptic boundary value problems that $\Lambda_\sigma$ is bounded. The Calder\'on problem is to determine $\sigma$ from the knowledge of $\Lambda_\sigma$.

As explained above, in this work we assume to be able to measure (an approximation of) the N-to-D map $\RR_\sigma: H_\diamond^{-1/2}(\partial\Om) \to H^{1/2} (\partial\Om)/ \C$ instead of the D-to-N map. If the current density on $\partial \Omega$ is set to $f \in H_\diamond^{-1/2}(\partial\Om)$, then the electrostatic potential $u\in H^1(\Om)/ \C$ uniquely solves the Neumann problem
\begin{eqnarray}
\label{eq:neumann}
  \left\{ \begin{array}{ll}
  \nabla \cdot (\sigma \nabla u) = 0  \quad & \rm{in} \ \Om \label{ody2},  \\[1mm]
  {\displaystyle \sigma\frac{\partial u}{\partial \nu}} = f  \quad &\rm{on} \ \partial \Om 
\end{array} \right.
\end{eqnarray}
interpreted in the weak sense.
Notice that \eqref{eq:neumann} does not fix the ground level of potential and thus the solution is unique only up to an additive constant,~i.e.,~up to the ground level of potential.
We then define
\begin{equation}\label{NDdef}
\RR_\sigma f = u|_{\partial\Om} \in H^{1/2}(\partial \Omega) / \C,
\end{equation}
where $u$ is the solution to \eqref{eq:neumann}.
It is easy to see that $\RR_\sigma$ is the inverse of $\Lambda_\sigma$, that is, 
\begin{equation}\label{iden}
\Lambda_\sigma \RR_\sigma = I \qquad {\rm and} \qquad \RR_\sigma \Lambda_\sigma  = I,
\end{equation}
where the identity map $I$ is that of $H_\diamond^{-1/2}(\partial\Om)$ in the first identity and that of $H^{1/2} (\partial\Om)/ \C$ in the second one. In particular, $\RR_\sigma$ is bounded.
Take note that the boundary quotient space $H^{1/2}(\partial\Om)/ \C$ can be replaced by its dual $H^{1/2}_{\diamond}(\partial\Om)$ at all occurrences in the above definitions if the ground level of potential is systematically chosen so that all electrostatic potentials have vanishing mean over $\partial \Omega$.

The direct reconstruction approach is to find equations that link $\Lambda_\sigma$ to $\sigma$ as well as to implement such a connection as a practical algorithm. The D-bar method originates from the work of Ablowitz, Nachman, Beals and Coifman \cite{Beals1981,Ablowitz1986} in a related inverse scattering problem. It is based on a \textit{non-linear Fourier transform}, where exponentially behaving \textit{complex geometric optics} (CGO) solutions of Faddeev \cite{Faddeev1966} are used. A \textit{boundary integral equation} deduced by R.~G.~Novikov \cite{Novikov1988} and Nachman \cite{Nachman1988a} links $\Lambda_\sigma$ to the CGO solutions. The two-dimensional Calder\'on problem was solved,~i.e.,~its uniqueness was deduced, using the D-bar method by Nachman \cite{Nachman1996} for conductivities in $W^{2,p}(\Om)$, $p>1$. The solution technique was later generalized by Brown and Uhlmann \cite{Brown1997} for less regular conductivities. Using D-bar techniques and quasiconformal maps, the two-dimensional Calder\'on problem was finally completely solved by Astala and P\ad iv\ad rinta \cite{Astala2006a} for $\sigma\in L^\infty_+(\Om)$.

The ideas of exploiting conformal transformations presented in this manuscript are arguably as useful for any of the aforementioned variants of the D-bar method and, more generally, for any reconstruction method where the D-to-N or N-to-D map is approximated  by a matrix in a truncated Fourier basis.
However, we will only numerically test these ideas using Nachman's original method.

\section{The D-bar method for EIT}\label{sec:Dbar}
The origin of the D-bar method lies with the direct solution to a certain inverse scattering problem formulated as the \textit{Gel'fand--Calder{\'o}n} (GC) problem \cite{Gelfand1961,Calder'on1980}: 
By writing $u = \sigma^{1/2}v$, \eqref{ody1} is transformed into a zero-energy Schr\od dinger equation
\begin{equation}\label{schrodinger}
(-\Delta+q)v = Ev \qquad \textrm{in }\R^2,\qquad q = \frac{\Delta\sqrt{\sigma}}{\sqrt{\sigma}},\quad E=0,
\end{equation}
which necessitates the smoothness assumption $\sigma\in W^{2,p}(\Om)$. Similarly, the equations governing {\em diffuse optical tomography} (DOT) and {\em acoustic tomography} (AT) can be recast in the form \eqref{schrodinger} but with negative and positive energies~\cite{deHoop2015, Tamminen2016}, respectively. 
We remark what the method of Astala and P\ad iv\ad rinta~\cite{Astala2006a} does not require transforming the conductivity equation into the form \eqref{schrodinger}, and in theory it thus  works with more general conductivities $\sigma\in L^\infty_+(\Omega)$. However, a number of numerical results indicate that Nachman's method remains functional even for discontinuous conductivities or if the smoothness assumptions are violated in some other way~\cite{Knudsen2007,Knudsen2008,ST2014}.

{\em Complex Geometrical Optics} (CGO) solutions $\psi$ are exponentially behaving solutions to the  Schr\od dinger equation \eqref{schrodinger}. They depend on the spatial variable $z\in\C$ as well as on a spectral variable $\zeta\in\C^2$ which satisfies $\zeta\cdot\zeta = E = 0$. Such $\zeta$ can be parametrized by a single complex parameter $k\in\C$. The exponential behavior is controlled by
\begin{equation}
\exp(-\rmi kz)\psi(z,k)-1\in W^{1,\tilde{p}}(\R^2), \qquad 1/\tilde{p} = 1/p-1/2.
\end{equation}
Nachman dubbed the potential $q$ appearing in \eqref{schrodinger} as of \textit{conductivity type} and proved that for such potentials the CGO solution $\psi(z,k)$ exists for all $k\in\C$. Thus for any $\sigma\in W^{2,p}(\Om)$ and $k\in\C$, one may define the scattering transform,~i.e.,~a certain \textit{nonlinear Fourier transform}, as
\begin{equation}\label{scat}
\T(k) = \int_{\Omega}e^{\rmi\kk{z}\kk{k}}q(z)\psi(z,k)\dee z.
\end{equation}
It can be shown that the function $\mu(z,k) = \psi(z,k)\exp(-\rmi kz)$ satisfies the following \textit{D-bar equation}:
\begin{equation}\label{Dbar}
\partial_\kk{k} \mu(z,k) = \frac{\T(k)e_{-k}(z)}{4\pi\kk{k}}\kk{\mu(z,k)},\qquad e_{-k}(z) = \exp(-\rmi(kz+\kk{k}\kk{z})),\quad k\neq 0,
\end{equation}
which is uniquely solvable in a certain weighted Sobolev space \cite{Nachman1996}. The conductivity can then be recovered via
\begin{equation}\label{sigmarecon}
\sigma(z) = \mu(z,0)^2.
\end{equation}
What remains to be explained is how the CGO solutions as well as the equations \eqref{scat} and \eqref{Dbar} are connected to the boundary measurement $\Lambda_\sigma$.

For simplicity, we assume in the following that $\sigma \equiv 1$ in some interior neighborhood of $\partial \Omega$. Let us define a single-layer operator 
\begin{equation}\label{Sk}
  (S_k \phi)(z) :=\int_{\DOm}G_k(z-y)\phi(y)\dee s(y),
\end{equation}
where $G_k$ is the \textit{Faddeev's Green function} 
\begin{equation}\label{Gk}
  G_k(z) := e^{\rmi kz}g_k(z),\qquad g_k(z) := \frac{1}{(2\pi)^2} \int_{\R^2}\frac{e^{\rmi z\cdot y}}{|y|^2+2k(y_1+\rmi y_2)}\dee y.
\end{equation}
For all $k\in\C\setminus\{0\}$, the trace of the CGO solution $\psi(\, \cdot \, ,k)|_{\DOm} \in H^{1/2}(\partial\Om)$ can be solved from the boundary integral equation
\begin{align}\label{BIE}
 \big(I + S_k(\Lambda_\sigma-\Lambda_1)\big)\psi(\,\cdot\,,k)|_{\DOm} = e^{{\rm i} k \, \cdot} ,
\end{align}
where $\Lambda_1$ is the D-to-N map corresponding to the homogeneous unit conductivity.

When solving \eqref{BIE}, the singularity in the kernel of the operator $S_k$ can be handled as follows; see~\cite{Siltanen1999} for more details. We introduce the harmonic function $H_k = G_k-G_0$, where $G_0(z) = -1/(2\pi)\log\abs{z}$ is the fundamental solution for the Laplacian, and note that $H_k(z) = H_1(kz)-1/(2\pi)\log\abs{k}$. This motivates defining yet another auxiliary function $\hat{H}_k(z) = H_1(kz)-H_1(0)$, for which $\hat{H}_k(0) = 0$. It follows that the kernel of $S_k$ can be restructured as
\begin{equation}
\label{eq:restruct}
G_k = (G_0 + \hat{H}_k) + \Big( H_1(0)- \frac{1}{2\pi} \log\abs{k} \Big).
\end{equation}
Since the latter term is independent of the spatial variable and it is known that $2 S_0 = \RR_1$, it is easy see that
$$
S_k =  \frac{1}{2} \RR_1 + \hat{\mathcal{H}}_k  \qquad {\rm on} \ H^{-1/2}_\diamond(\partial \Omega),
$$
with $\hat{\mathcal{H}}_k$ being the integral operator defined by the well-behaving kernel $\hat{H}_k$. In particular,
\begin{align}
I + S_k(\Lambda_\sigma-\Lambda_1) &= I + \Big(\frac{1}{2} \RR_1 + \hat{\mathcal{H}}_k \Big) (\Lambda_\sigma-\Lambda_1) \nonumber \\[1mm]
&= \frac{1}{2}(I + \RR_1\Lambda_\sigma )+ \hat{\mathcal{H}}_k(\Lambda_\sigma-\Lambda_1). \label{BIEoperator}
\end{align}
The operator $\hat{\mathcal{H}}_k$ includes the Faddeev's Green function $G_k$ as a part of its kernel.

The scattering transform \eqref{scat} also admits an alternative form
\begin{align}\label{scatDN}
  \T(k)=\int_{\partial\Om}e^{\rmi\bar{k}\bar{z}} (\Lambda_\sigma-\Lambda_1)\psi(z,k)\dee s(z),
\end{align}
which can be evaluated after solving $\psi(z,k)|_{\DOm}$ from \eqref{BIE}. 
Subsequently, the D-bar equation \eqref{Dbar} can be written as a Lippman--Schwinger type integral equation, uniquely solvable in $k\in\C$ for each fixed $z\in \Om$:
\begin{align}\label{IE}
\mu(z,k)=  1+\frac{1}{(2\pi)^2}\int_{\R^2}\frac{\T(k^\prime)}{(k-k^\prime)\bar{k}^\prime} e^{\rmi(k^\prime z+\kk{k^\prime}\kk{z})}\kk{\mu(z,k^\prime)}\dee k^\prime_1 \dee k^\prime_2.
\end{align}
Finally, \eqref{sigmarecon} provides the solution to the two-dimensional Calder\'on problem.

\subsection{Standard numerical implementation of the D-bar method}
The D-bar method has a standard numerical implementation; see \cite{Siltanen2001,Mueller2002,Knudsen2004} and \cite{Isaacson2004,Isaacson2006} for its applications to simulated and experimental data, respectively.

Let the Lipschitz continuous, $|\DOm|$-periodic function $\gamma: \R \to \C$ define a parametrization of $\partial \Omega$ with respect to its arclength, that is, $|\gamma'| \equiv 1$ almost everywhere and
\begin{equation}
\DOm = \{z\in\C \ | \ z = \gamma(s), \ s \in \R\}.
\end{equation}
The corresponding $L^2(\partial \Omega)$-orthonormal \textit{truncated Fourier basis} is defined via
\begin{equation}
\label{Fbasis}
\phi_n(s) = \frac{1}{\sqrt{\abs{\DOm}}}\exp\Big(\rmi \frac{2\pi n}{\abs{\DOm}} s \Big), \qquad n=-N,...,N,
\end{equation}
as functions of the arclength parameter. The linear operators $\RR_\sigma,\RR_1,\Lambda_\sigma,\Lambda_1$ and $\hat{\mathcal{H}}_k$, needed in \eqref{BIEoperator} and \eqref{scatDN}, are approximated by matrices $\mtx{R}_\sigma, \mtx{R}_1,\mtx{L}_\sigma,\mtx{L}_1$ and $\hat{\mtx{H}}_k$, respectively, with respect to the basis \eqref{Fbasis}. To be more precise,
a bounded linear operator (see,~e.g.,~\cite[p.~20]{Grisvard1985})
$$
A: \ H^{s}(\DOm)\to H^{r}(\DOm), \qquad s,r \in (-1,1),
$$
is replaced in the numerical algorithm by the matrix $\mtx{A} \in \C^{(2N+1) \times (2N+1)}$ defined through
\begin{equation}\label{Amn}
\mtx{A}_{m,n} = \langle A\phi_n, \overline{\phi}_m\rangle_{\partial \Omega}, \qquad m,n=-N,\dots,N,
\end{equation}
where $m,n=-N$ refers to the first row/column of the matrix. If $\vct{f}_n = \langle f, \overline{\phi}_n \rangle_{\partial \Omega}$, $n = -N, \dots, N$, are the Fourier coefficients of an input $f \in H^{s}(\DOm)$, then the Fourier coefficients of $Af$ can be approximated as
$$
\langle A f, \overline{\phi}_m \rangle  \approx [\mtx{A}\vct{f}]_m, \qquad m=-N,\dots, N.
$$ 
It is straightforward to deduce that the appropriate interpretation of $\mtx{A}$ as a finite-dimensional mapping from $H^{s}(\DOm)$  to $H^{r}(\DOm)$ converges to $A:  H^{s}(\DOm)\to H^{r}(\DOm)$ pointwise as $N$ goes to infinity (see,~e.g.,~\cite{Saranen2002}). Furthermore, the convergence occurs in the operator norm if $A: H^{s}(\DOm)\to H^{r}(\DOm)$ is compact. The case of mean-free spaces, i.e., $A: H^{s}(\DOm)/\C \simeq H^{s}_\diamond(\DOm)\to H^{r}_\diamond(\DOm) \simeq H^{r}(\DOm)/\C$ can be treated by simply dropping out the constant basis function corresponding to $n=0$. 

To simulate data in a bounded domain,\footnote{See Appendix~\ref{sec:halfplane} for the case of the upper half-plane.} we compute the necessary direct solutions $u$ of \eqref{ody2} in order to construct the N-to-D matrix $\mtx{R}_\sigma$ by a {\em finite element method} (FEM) with piecewise linear basis functions. The corresponding D-to-N matrix $\mtx{L}_\sigma$ can be obtained by inverting $\mtx{R}_\sigma$, though this process becomes unstable for large $N$ since $\RR_\sigma$ is compact and $\Lambda_\sigma$ unbounded on $L^2_\diamond(\partial \Omega)$. 
To compute the matrix $\hat{\mtx{H}}_k$ approximating the operator $\hat{\mathcal{H}}_k$ in \eqref{BIEoperator} we need an algorithm for evaluating $G_k(z)$ of \eqref{Gk}. In the zero energy case this is trivial, as we can use the property $g_k(z) = g_1(kz)$ and \cite[(3.10)]{Boiti1987} to compute $g_1$ using MATLAB's built-in exponential-integral function {\tt expint}.

Solving \eqref{BIE} could be skipped by replacing $\psi(z,k)$ with $\exp(\rmi kz)$ in \eqref{scatDN} without considerably affecting the reconstructions (see,~e.g.,~\cite{Siltanen2000}) because the two functions have the same asymptotic behavior as  $\abs{z}$ tends to infinity. However, here we solve the full problem by approximating \eqref{BIE} as a matrix equation (cf.~\eqref{BIEoperator})
\begin{equation}\label{BIE2}
\Big (\frac{1}{2}(\mtx{I}+\mtx{R}_1\mtx{L}_\sigma)+\hat{\mtx{H}}_k(\mtx{L}_\sigma-\mtx{L}_1) \Big)\vct{p}_k = \vct{e}_k
\end{equation} 
in the truncated Fourier basis.
Subsequently, an approximate scattering transform is computed as
\begin{equation}\label{scatDN2}
\vct{t}(k) = \int_{\partial\Om}e^{\rmi\bar{k}\bar{z}}
    \mathcal{F}^{-1}\big((\mtx{L}_\sigma-\mtx{L}_1)\vct{p}_k\big)(z) \, \dee s(z),
\end{equation}
where $\mathcal{F}^{-1}$ corresponds to evaluating the function with the Fourier coefficients $(\mtx{L}_\sigma-\mtx{L}_1)\vct{p}_k$. The integral in \eqref{scatDN2} is approximated by the trapezoidal rule with respect to a dense enough discretization of the arclength parameter $s\in[0,\abs{\DOm})$.

In practice, the discretized boundary integral equation \eqref{BIE2} is not (stably) solvable for all $k \in \C$; typically, problems arise when $|k|$ is large and/or when the data contain noise. 
However, the D-bar method can be regularized by truncating the scattering transform, that is, introducing
\begin{align}\label{scatR}
    \T_R(k)=
  \begin{cases}
\T(k)& \mbox{for }|k|<R, \\ 
0 & \text{otherwise},
  \end{cases}
\end{align}
with $R>0$ suitably defined as a decreasing function of the noise level~\cite{Knudsen2009}. In this work, we resort to a somewhat more general truncation
\begin{align}\label{scatRc}
    \T_{R,c}(k)=
  \begin{cases}
\T_c(k)& \mbox{for }|k|<R,\\ 
0 & \text{otherwise},
  \end{cases}
\end{align}
where $\T_c(k) = 0$ when $\abs{\T(k)}>c$ with $c>0$ being a suitable cut-off parameter. We have $\abs{\T_{R,c}(k)}\leq \abs{\T_{R}(k)}$ for every $k \in \C$ so the convergence results of \cite{Knudsen2009} still hold. However, the analysis on why \eqref{scatRc} would be a better truncation than \eqref{scatR} is nontrivial and absent to the authors' best knowledge. In our numerical experiments, the values of $R, c >0$ are chosen after visually investigating the right-hand side of \eqref{scatDN2} on an extensive grid over $k$ in the complex plane.

The integral equation \eqref{IE} can be tackled,~e.g.,~by a modified Lippman--Schwinger solver originating from Vainikko's work~\cite{Vainikko2000}: The equation is periodized and solved using GMRES as explained in \cite{Knudsen2004}. The resulting solution is (an approximation of) the function $\mu(z,\,\cdot)$ at the investigated reconstruction point $z \in \Omega$ and a grid of values over the second variable,~i.e.,~$k$. We then set $k=0$ and get the reconstructed conductivity at $z$ via~\eqref{sigmarecon}.

\subsection{The case of unit disk}
\label{sec:recon}

If the D-bar method is implemented in a domain that is not a disk, the reference N-to-D map $\mathcal{R}_1$ (or $\Lambda_1$) must also be approximated numerically by,~e.g.,~FEM; cf.~\eqref{BIE2} and \eqref{scatDN2}. One benefit in transforming the reconstruction problem to the unit disk is that there these auxiliary maps admit trivial diagonal representations, which makes the implementation of the D-bar algorithm more efficient.

Indeed, the reference D-to-N map $\tilde{\Lambda}_1$ corresponding to the unit disk $D := \{z\in\C: \abs{z}<1\}$ with homogeneous unit conductivity can be characterized in terms of the Fourier basis \eqref{Fbasis} as
\begin{equation}
\label{eq:Lunitdisk}
\tilde{\Lambda}_1\tilde{\phi}_n = \abs{n}\tilde{\phi}_n,\qquad \abs{n}>1.
\end{equation}
Here and in what follows, we denote a Fourier basis function by $\tilde{\phi}_n$ if it corresponds to $\Omega = D$, that is, $\tilde{\phi}_n(\theta) = 1/\sqrt{2\pi}\exp(\rmi n\theta)$, with $\theta$ being the polar angle. Such `tilde notation' is systematically adopted to indicate which entities are related to $D$. As a matrix with respect to truncated Fourier basis, the relations \eqref{eq:Lunitdisk} takes the form 
\begin{equation}
\label{eq:Ldisk}
\tilde{\mtx{L}}_1 = \textrm{diag}\big(N,N-1,\ldots,2,1,1,2,\ldots,N-1,N\big) \in \R^{2N \times 2N},
\end{equation}
where the constant basis function has been dropped. This obviously also induces a diagonal representation for the reference N-to-D matrix $\tilde{\mtx{R}}_1$.

The following lemma lists other symmetries that streamline the implementation of the D-bar method in the unit disk by reducing the computations required for forming the matrix $\hat{\mtx{H}}_k$ appearing in \eqref{BIE2}. These results are formulated for the matrix $\mtx{S}_k$ approximating the single-layer operator \eqref{Sk}, but they also hold for $\hat{\mtx{H}}_k$ due to the relation
$$
\hat{\mtx{H}}_k = \mtx{S}_k-\frac{1}{2}\mtx{R}_1+\frac{1}{2\pi}\log(\abs{k}) \, \mtx{I}
$$
induced by \eqref{eq:restruct}.
The proof is omitted since the claims are easy consequences of the symmetries exhibited by the Faddeev Green's function $G_k$.

\begin{Lemma}
For all $k\in\C\setminus\{0\}$,
\begin{equation*}
\mtx{S}_{-k} = \kk{\mtx{S}_k^{{\rm T}}}.
\end{equation*}
If $\Omega$ is a disk, then also 
\begin{align*}
\mtx{S}_{\kk{k}} = \mtx{S}_{k}^{\rm T} \qquad {\rm and} \qquad
\mtx{S}_{k} = \mtx{E}_\alpha \odot \mtx{S}_{\abs{k}},
\end{align*}
where $\odot$ denotes elementwise multiplication of matrices, $\alpha$ is the polar angle of $k\in\C$ and $(\mtx{E}_\alpha)_{n,m} = \exp(\rmi\alpha(m-n))$.
\end{Lemma}

\section{Conformal transformation of the domain}
\label{sec:conformal}
Let $D = \{z\in\C: \abs{z}<1\}$ still be the open unit disk. According to the Riemann mapping theorem, there is a biholomorphic/conformal map $\Phi: \Omega\to D$. Since $\partial \Omega$ is Lipschitz, $\Phi$ extends to a homeomorphism $\Phi: \overline{\Omega}\to \overline{D}$~\cite{Pommerenke92}. The inverse of $\Phi$ is denoted by $\Psi$. 
We denote by $J_\Phi$ the Jacobian of $\Phi$ (interpreted as a map from $\R^2$ to itself)  and by $|\Phi'| = \sqrt{\det J_\Phi}$ the absolute value of its (complex) derivative.
The following lemma, the basic principle of which is well-known, relates the solution of \eqref{eq:neumann} to that of
\begin{eqnarray}
\label{eq:neumann2}
  \left\{ \begin{array}{ll}
  \nabla \cdot (\tilde{\sigma} \nabla \tilde{u}) = 0  \quad & {\rm in} \ D,  \\[1mm]
  {\displaystyle \tilde{\sigma} \frac{\partial \tilde{u}}{\partial \nu}} = \tilde{f}  \quad &{\rm on} \ \partial D
\end{array} \right.
\end{eqnarray}
for a certain $\tilde{f} \in H^{-1/2}_\diamond(\partial D)$ and $\tilde{\sigma} := \sigma \circ \Psi$.

\begin{Lemma}
\label{changevar}
Let $u \in H^1(\Omega) / \C$ be the solution of \eqref{eq:neumann} for some $f \in H^{-1/2}_\diamond(\partial \Omega)$. Then, $\tilde{u} := u \circ \Psi \in H^1(D) / \C$ is the solution of~\eqref{eq:neumann2} for $\tilde{f} \in H^{-1/2}_\diamond(\partial D)$ defined by
\begin{equation}
\label{eq:pushforward}
\langle \tilde{f},  \tilde{g} \rangle_{\partial D} \, = \, \langle f, \tilde{g} \circ \Phi|_{\partial \Omega} \rangle_{\partial \Omega}
\end{equation}
for all $\tilde{g} \in H^{1/2}(D)/\C$.
\end{Lemma}

\begin{proof}
We mimic the proof of \cite[Lemma~A.3]{Seiskari14} and start by showing the linear map $L_\Phi: \tilde{v} \mapsto \tilde{v}\circ \Phi$ is bounded from $H^1(D) / \C$ to $H^1(\Omega) / \C$. Let us denote by $D_\epsilon$ the open, origin-centered disk of radius $1-\epsilon$ and by $\chi_{D_\epsilon}$ its characteristic function. For any $\tilde{v} \in H^1(D) / \C$, it holds that
\begin{equation}
\begin{aligned}
\int_{D} \chi_{D_\epsilon} | \nabla \tilde{v} |^2 \, {\rm d} z  
&= \int_{D_\epsilon} \nabla \tilde{v} \cdot \nabla \overline{\tilde{v}}  \, {\rm d} z  \\ 
&= \int_{\Psi(D_\epsilon)} J_\Phi^{\rm -T} \nabla (\tilde{v}\circ\Phi) \cdot J_\Phi^{\rm -T} \nabla (\overline{\tilde{v}\circ\Phi}) \, \abs{\Phi'}^2 \, {\rm d}w  \\
&= \int_{\Psi(D_\epsilon)} \nabla (\tilde{v}\circ\Phi)^{\rm T} J_\Phi^{\rm -1}  J_\Phi^{\rm -T} \nabla (\overline{\tilde{v}\circ\Phi}) \, \abs{\Phi'}^2 \, {\rm d}w \\
&= \int_{\Omega} \chi_{\Psi(D_\epsilon)}|\nabla (\tilde{v}\circ\Phi)|^2  \, {\rm d}w
\end{aligned}\label{jee}
\end{equation}
where we used $J_\Phi^{\rm -1} J_\Phi^{\rm -T} = (1/\abs{\Phi'}^2) I$, a consequence of the Cauchy--Riemann equations. Since the integrands of the first and the last term of \eqref{jee} are clearly increasing as $\epsilon \to 0$, the Lebesgue monotone convergence theorem gives 
$$
\| \nabla(\tilde{v} \circ \Phi) \|_{L^2(\Omega)} = \| \nabla\tilde{v}\|_{L^2(D)},
$$
which shows the boundedness of $L_\Phi: H^1(D) / \C \to H^1(\Omega) / \C$ by virtue of the Poincar\'e inequality.  (Note that $L_\Phi$ is actually a homeomorphism as its inverse is obviously defined by $v \mapsto v\circ \Psi$.) 

As the trace map is continuous from $H^1(\Omega)/\C$ to $H^{1/2}(\partial \Omega)/\C$~\cite{Hyvonen2004} and $\Phi|_{\partial \Omega}: \partial \Omega \to \partial D$ is a homeomorphism~\cite{Pommerenke92}, it follows,~e.g.,~by the density of the embedding $C^\infty(\overline{D}) \hookrightarrow H^1(D)$~\cite{Grisvard1985}, that any $\tilde{v} \in H^1(D)/\C$ satisfies
\begin{equation}
\label{eq:reuna}
(\tilde{v} \circ \Phi)|_{\partial \Omega} = \tilde{v}|_{\partial D} \circ \Phi|_{\partial \Omega}
\end{equation}
modulo an additive constant almost everywhere on $\partial \Omega$.
Let $\tilde{\gamma}^{-1}: H^{1/2}(\partial D)/\C \to H^1(D)/\C$ be a bounded right-inverse for the (quotient) trace operator in $D$~(cf.~\cite{Grisvard1985}). By \eqref{eq:reuna}, the boundedness of the quotient trace map and the continuity of $L_\Phi$, an arbitrary $\tilde{g} \in H^{1/2}(\partial D)/\C$ satisfies
\begin{align}
\label{eq:reunahom}
\| \tilde{g} \circ \Phi|_{\partial \Omega} \|_{H^{1/2}(\partial \Omega)/\C} &= 
\big\| \big( (\tilde{\gamma}^{-1} \tilde{g}) \circ \Phi \big) \big|_{\partial \Omega} \big\|_{H^{1/2}(\partial \Omega)/\C} \leq C\| (\tilde{\gamma}^{-1} \tilde{g}) \circ \Phi \|_{H^1(\Omega) /\C}  \nonumber \\[1mm]
&\leq C \| \tilde{\gamma}^{-1} \tilde{g} \|_{H^1(D) /\C} \leq C \| \tilde{g} \|_{H^{1/2}(\partial D)/\C}. 
\end{align}
In particular,
$$
| \langle \tilde{f}, \tilde{g} \rangle_{\partial D} | \leq 
C \| f \|_{H^{-1/2}(\partial \Omega)} \| \tilde{g} \|_{H^{1/2}(\partial D) /\C},
$$
where the fact $\langle \tilde{f}, 1 \rangle_{\partial D} = 0$ was implicitly used.
Together with the obvious linearity of $\tilde{f}$, this demonstrates that $\tilde{f} \in H^{-1/2}_\diamond(\partial D)$.

A similar computation as \eqref{jee},  H\od lder's inequality and the Lebesgue dominated convergence theorem can be used to show that
$$
\int_{D} \tilde{\sigma} \nabla \tilde{u} \cdot \nabla \overline{\tilde{v}} \, {\rm d} z \,= \, \int_{\Omega} \sigma \nabla u \cdot \nabla ( \overline{\tilde{v}\circ\Phi}) \, {\rm d} w \, = \, \langle f, (\overline{\tilde{v}\circ\Phi})|_{\partial \Omega} \rangle_{\partial \Omega} \quad {\rm for} \ {\rm all} \ \tilde{v} \in H^1(D)/\C,
$$
where the second equality corresponds to the variational formulation of \eqref{eq:neumann} with the test function $\tilde{v} \circ \Phi \in H^1(\Omega)/\C$~\cite{Grisvard1985}. It thus follows from \eqref{eq:reuna} and the definition of $\tilde{f}$ that actually
$$
\int_{D} \tilde{\sigma} \nabla \tilde{u} \cdot \nabla \overline{\tilde{v}} \, {\rm d} z = \langle \tilde{f}, \overline{\tilde{v}}|_{\partial D}  \rangle_{\partial D} \quad {\rm for} \ {\rm all} \ \tilde{v} \in H^1(D)/\C,
$$
which is the variational formulation of \eqref{eq:neumann2} and completes the proof. 
\end{proof}

\begin{Remark}
If $f \in L^2_\diamond(\partial \Omega)$ and the boundary $\partial \Omega$ is regular enough, one can make the interpretation
$$
\langle f, \tilde{g}\circ \Phi|_{\partial \Omega} \rangle_{\partial \Omega} = 
\int_{\partial \Omega} f (\tilde{g} \circ \Phi) \, {\rm d} s 
= \int_{\partial D} |\Psi'| (f\circ \Psi) \, \tilde{g}  \, {\rm d} s
$$
in \eqref{eq:pushforward}. This holds,~e.g.,~if the Lipschitz boundary $\partial \Omega$ consists of a finite number of $C^{1,\alpha}$ smooth arcs, with $0 < \alpha \leq 1$, because then $\Psi$ has a continuously differentiable extension onto the open $C^{1,\alpha}$ arcs~\cite{Pommerenke92}, allowing a change of variables one arc at a time. In other words, 
\begin{equation}
\label{eq:piecewsmooth}
\tilde{f} \, = \, |\Psi'| (f\circ \Psi) \qquad {\rm or} \qquad 
f = |\Phi'| (\tilde{f} \circ \Phi)
\end{equation}
in the sense of $L^2(\partial D)$ and $L^2(\partial \Omega)$, respectively.
\end{Remark}

\begin{figure}
\begin{picture}(120,120)
	\put(10,50){\includegraphics*[width=4cm]{\imagepath 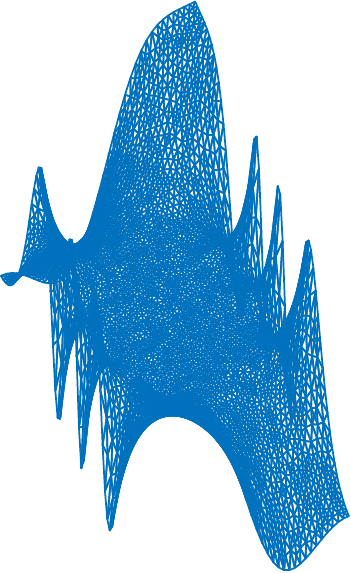}}
\put(60,50){\includegraphics*[width=4cm]{\imagepath 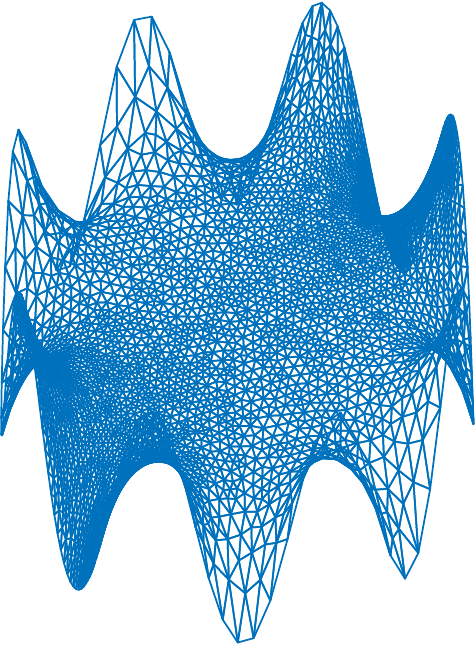}}
\put(10,0){\includegraphics*[width=4cm]{\imagepath 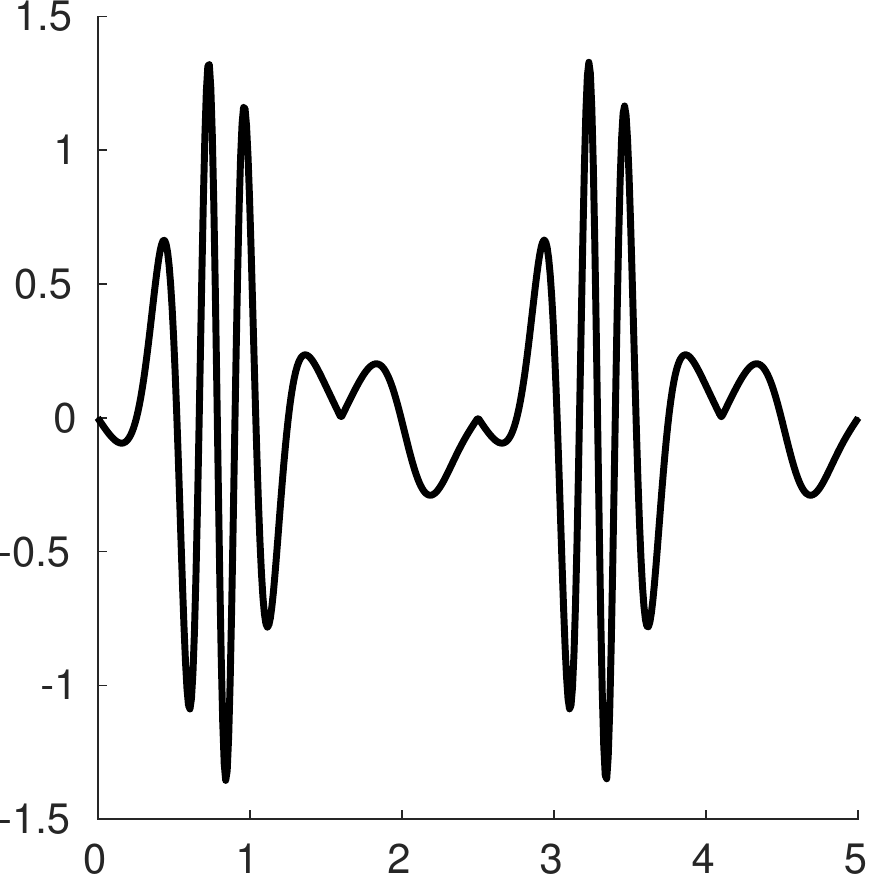}}
\put(60,0){\includegraphics*[width=4cm]{\imagepath 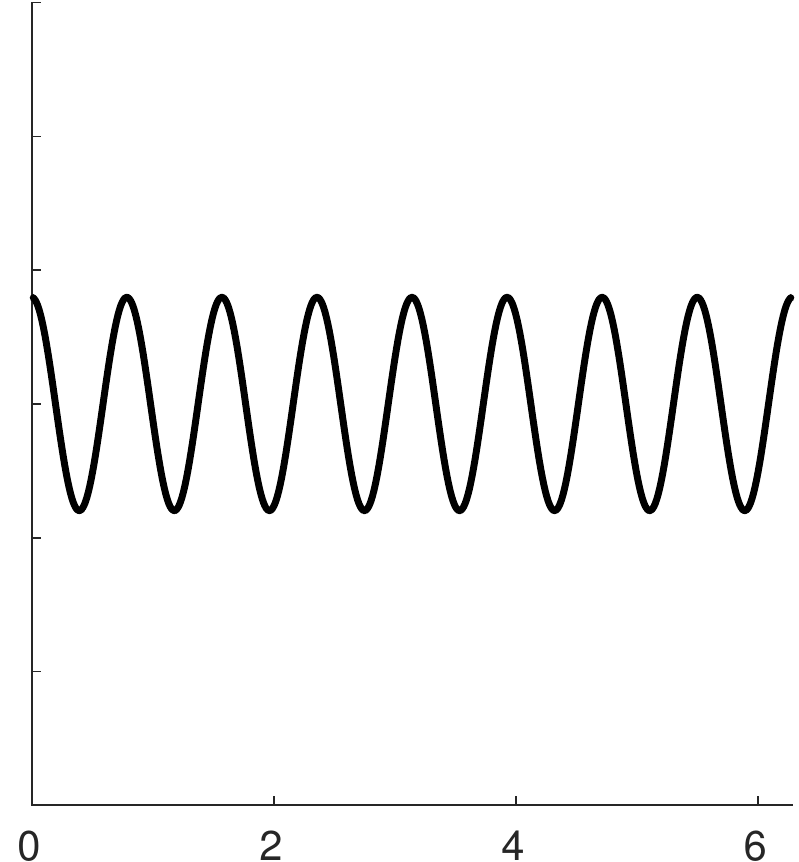}}

\end{picture}
\caption{Top left: The real part of a FEM approximation for the solution $u$ to \eqref{ody2} in a homogeneous rectangular domain $\Omega$ with the boundary current density $\varphi_8$ of \eqref{eq:trans} corresponding to a Schwarz--Christoffel map $\Phi: \Omega \to D$. Top right: The real part of $\tilde{u} = u \circ \Psi$ in $D$. Bottom left: The real part of the true current density $\varphi_8$ as a function of the arclength parameter on $\partial \Omega$. Bottom right: The real part of the virtual current density $\tilde{\phi}_8$ as a function of the polar angle on $\partial D$. \label{fig:transformation}}
\end{figure}

As we want to use the truncated Fourier current basis $\tilde{f} = \tilde{\phi}_n$, $1\leq |n| \leq N$, on the boundary of the virtual reconstruction domain $D$, the transformations \eqref{eq:piecewsmooth} indicates that one should drive the (true) current patterns
\begin{equation}
\label{eq:trans}
\varphi_n \, = \, |\Phi'| (\tilde{\phi}_n\circ \Phi), \qquad 1\leq |n| \leq N,
\end{equation}
through $\partial \Omega$. These are called {\em conformally
transformed Fourier basis} functions and were mentioned in \cite{Garde2016}.
See Figure \ref{fig:transformation} for an example of the boundary current transformation \eqref{eq:trans} in a rectangular $\Omega$. 

In the following, $\Phi$ is often a composition of two conformal maps: one being a Schwarz--Christoffel map produced by the SC-Toolbox~\cite{Driscoll1996} sending a given polygon onto the unit disk, the other being a M\"obius transform magnifying a certain ROI within the unit disk. (Take note that such a composition is also itself a Schwarz--Christoffel map.)

\subsection{M\od bius transforms}
\label{sec:mobius}
\begin{figure}
\begin{picture}(120,56)
\put(0,0){\includegraphics*[width=5cm]{\imagepath 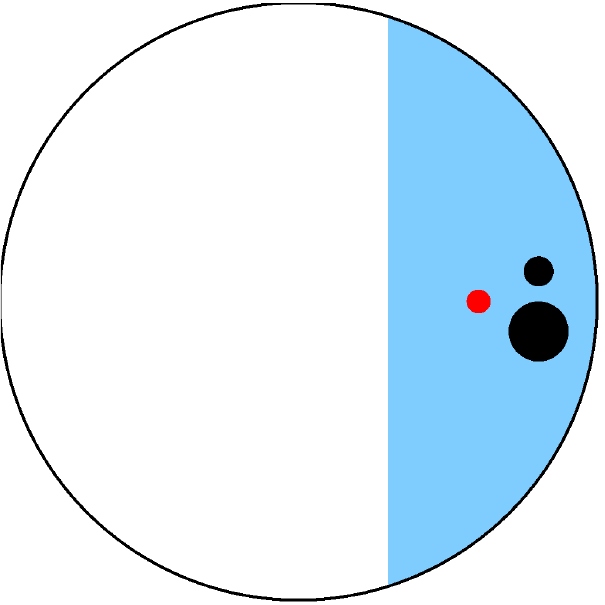}}
\put(60,0){\includegraphics*[width=5cm]{\imagepath 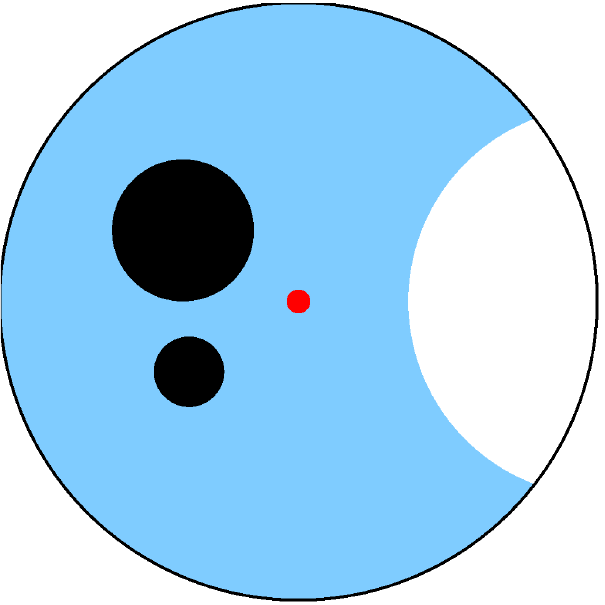}}
\put(42,48){$\sigma(z)\quad \longrightarrow \quad \tilde{\sigma}(z)$}
\put(52,43){$\MM_a$}
\end{picture}
\caption{Left: The true conductivity, with the blue area representing the ROI that contains two black anomalies. The red dot indicates the M\od bius parameter $a=0.6$. Right: The conformally mapped conductivity $\tilde{\sigma} = \sigma \circ \MM_a^{-1}$, with the image of the ROI under $\MM_a$ presented by blue color. 
\label{fig:mobius}}
\end{figure}

All conformal maps of $D$ onto itself are given by the M\od bius transforms of the form
\begin{equation}\label{mobius}
\MM_a(z) = e^{\rmi \alpha}\frac{z-a}{\kk{a}z-1},
\end{equation}
where $\alpha\in\R$ determines a rotation and $a\in D$ is a complex number that is mapped to the origin. In the following, we systematically choose $\alpha = 0$. Such M\od bius transforms can be used to magnify a certain part of the unit disk as illustrated by Figure~\ref{fig:mobius} for $a=0.6$. 
With such a choice, the inhomogeneities inside the ROI close to the right edge of $D$ are magnified. 

Our leading idea is that the virtual conductivity on the right in Figure~\ref{fig:mobius} is easier to reconstruct than the true one on the left if employing a {\em finite number of standard} Fourier basis functions,~i.e.,~the Fourier basis with respect to the polar angle $\{\tilde{\phi}\}_{1\leq |n|\leq N}$, as the input current densities on the respective boundaries. In other words, if it is a priori known that the ROI lies close to the right-hand edge of the true domain, it is beneficial to employ the {\em conformally transformed  Fourier basis functions} \eqref{eq:trans}, with $\Phi = \MM_a$, as the true input current densities for the true domain on the left, which corresponds to the use of the standard Fourier current basis for the virtual domain on the right. The standard, efficient numerical implementation of the D-bar algorithm in the unit disk described in Section~\ref{sec:recon} can then be applied to the virtual data on the boundary of the virtual domain (cf.~Lemma~\ref{changevar}), and the obtained reconstruction can be subsequently mapped by $\Psi$ back to the true domain. Our numerical experiments in Section~\ref{sec:numerics} demonstrate this procedure results in a more accurate reconstruction in the ROI. (Take note that the true domain does not need to be the unit disk nor the conformal map a M\"obius transform when magnifying a ROI.)

In Appendix~\ref{sec:point}, it is shown that relative continuum measurements for the true domain with the conformally transformed Fourier basis currents \eqref{eq:trans} can be approximated accurately by using a finite number of pointlike electrodes (cf.~\cite{Hanke2011}) if their locations are chosen to be the preimages under $\Phi$ of a uniform grid on the boundary of the virtual domain. In other words, the use of a conformally transformed Fourier basis  corresponds in practice to employing a certain nonuniform electrode grid that is densest close to the ROI.

\section{Numerical tests}
\label{sec:numerics}
Our numerical experiments are divided into three sections. In the first part, we investigate how reconstructions of a certain conductivity phantom in polygonal domains depend on whether the D-bar method is implemented directly in the true domain by resorting to the standard Fourier current basis \eqref{Fbasis} or in virtual domain (the unit disk) by using the conformally transformed Fourier current basis \eqref{eq:trans} on the boundary of the true domain. These tests do not consider the magnification of a ROI, but only aim at demonstrating that both approaches produce approximately as good reconstructions. This suggests that the choice between the two options should be made by comparing their computational expense.

In the second part, we numerically demonstrate that it is possible to produce more accurate reconstructions of a given ROI by magnifying it with a suitable conformal map. We consider two measurement geometries: In the first one, $\Omega = D$, meaning that the magnification is performed by a mere M\"obius transform. In the second case, $\Omega$ is a polygon and the needed conformal transform $\Phi$ is thus a Schwarz--Christoffel map \cite{Trefethen2002} (or a composition of a Schwarz--Christoffel map and a M\"obius transform as in our numerical implementation).

In the third part, we consider EIT measurements in a half-plane setting; simulation of (relative) EIT data for the half-plane is considered in Appendix~\ref{sec:halfplane}. The idea is to map the upper half-plane onto the unit disk by utilizing a family of M\"obius transforms, each of which magnifies different details in the target conductivity. By forming the corresponding reconstructions in the unit disk and mapping them back to the half-plane, one obtains reconstructions that highlight different details in the conductivity of the upper half-plane. In practice, the use of such a family of M\"obius transforms can be interpreted as moving an array of electrodes along the horizontal axis; see Appendices~\ref{sec:point} and \ref{sec:halfplane}.

All measurement data for bounded target domains are simulated by solving the Neumann problem \eqref{eq:neumann} for the employed current patterns by FEM with piecewise linear basis functions. For the half-plane geometry, the (relative) EIT data (for the virtual domain) are simulated by resorting to layer potential techniques, as explained in Appendix~\ref{sec:halfplane}. In both cases, the utilized discretizations are so dense that the (relative) amount of numerical noise is assumed to be insignificant for the employed lowish spatial Fourier frequencies; if not stated otherwise, we use the truncation index $N=16$, for all families of Fourier-like current patterns (cf.~\eqref{Fbasis} and \eqref{eq:trans}).

For each considered reconstruction domain, the $\tilde{\mtx{H}}_k$ matrix needed in \eqref{BIE2} is formed on a uniform grid of $128 \times 128$ points over the square $-12\leq k_1,k_2\leq12$, $k=k_1+\rmi k_2$. To evaluate \eqref{scatDN2}, we use the trapezoidal rule with 256 uniform nodes. The integral equation \eqref{IE} is discretized over a uniform $k$-grid of $2^7\times 2^7$ points. Our choice for the truncation parameters $R$ and $c$ in \eqref{scatRc} is somewhat heuristic: We first plot the scattering transform over the full $k$-grid mentioned above, then pick a reasonable $c>0$, and finally choose the smallest $R \in \N$ for which the connected area where $\T_{R,c}= \T$ around the origin is as large as possible; see Figure~\ref{fig:scats} for an example.

Finally, note that all considered conductivities are piecewise constant, which contradicts the $W^{2,p}$ assumption of Nachman's method \cite{Nachman1996} but anyway results in reasonable reconstructions (cf.~\cite{Knudsen2007,Knudsen2008,ST2014}).

\begin{figure}
\begin{picture}(120,85)
\epsfxsize=12cm
\put(0,40){\includegraphics*[width=3.8cm]{\imagepath 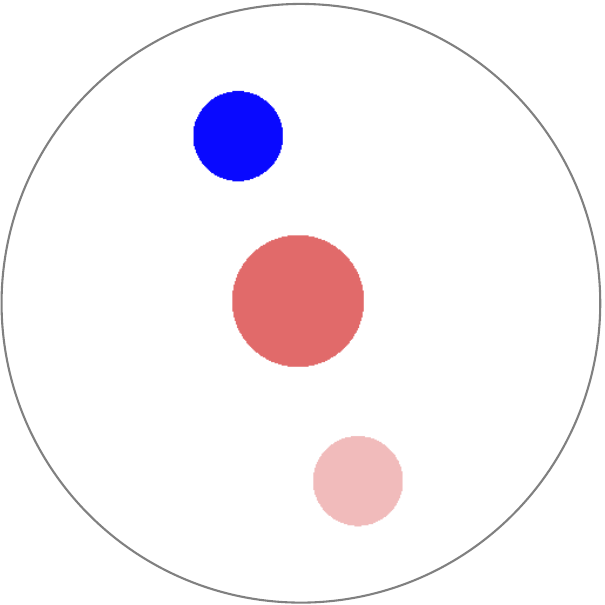}}
\put(40,40){\includegraphics*[width=3.8cm]{\imagepath 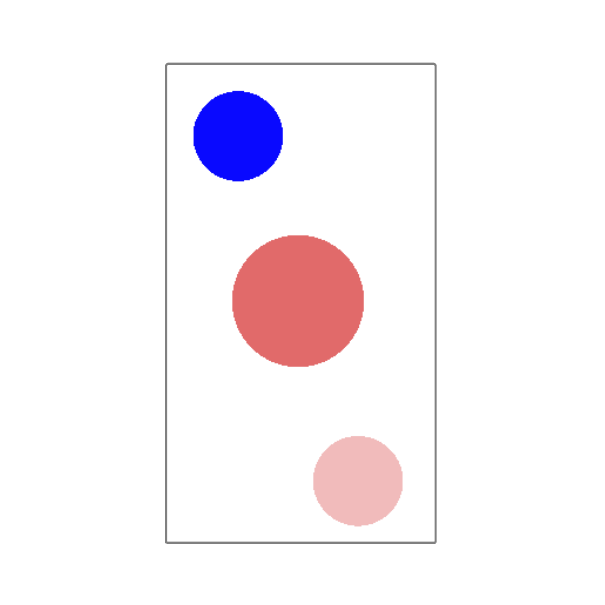}}
\put(80,40){\includegraphics*[width=3.8cm]{\imagepath 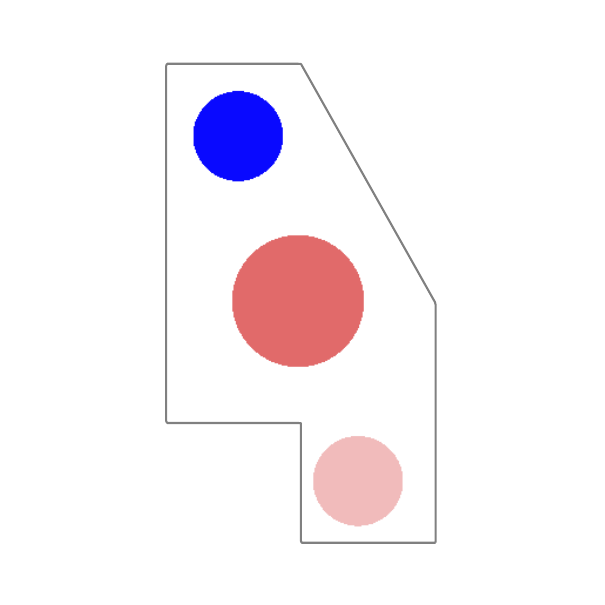}}

\put(40,0){\includegraphics*[width=3.8cm]{\imagepath 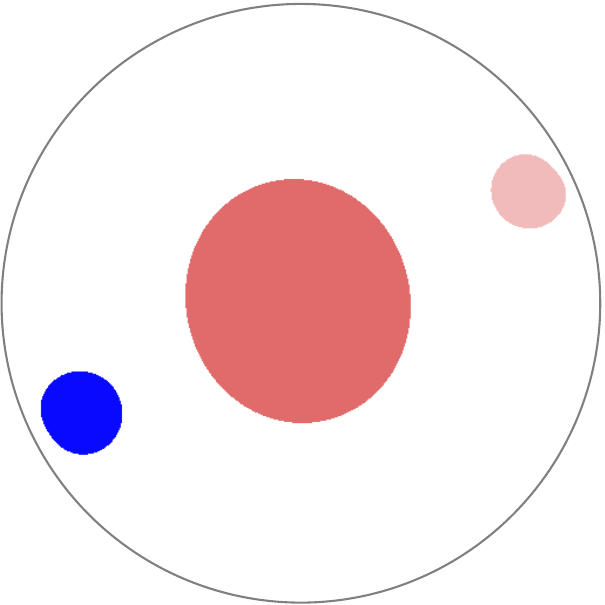}}
\put(80,0){\includegraphics*[width=3.8cm]{\imagepath 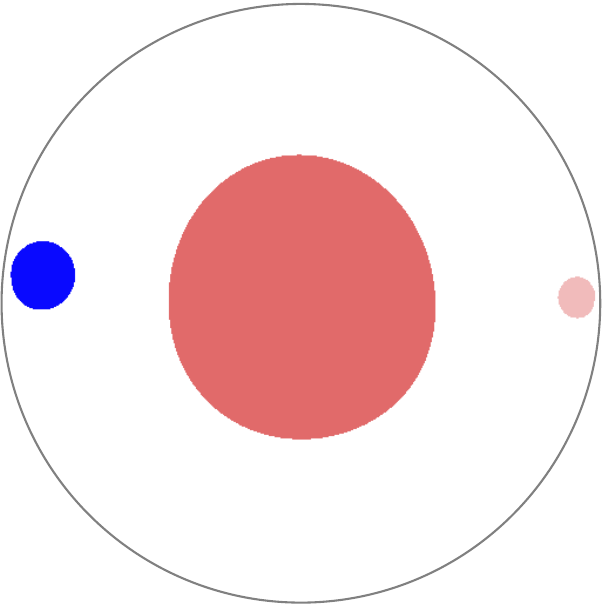}}

\put(8,80){$\Omega_1$, unit disk}
\put(48,80){$\Omega_2$, rectangle}
\put(88,80){$\Omega_3$, polygon}
\put(74,38){$\tilde{\sigma}(z)$}
\end{picture}
\caption{Top row: the target conductivity in three different domains. Bottom row: the virtual conducitivities,~i.e.,~the target conductivities transformed to the unit disk by a Schwarz--Christoffel map that fixes the origin. \label{fig:domains}}
\end{figure}
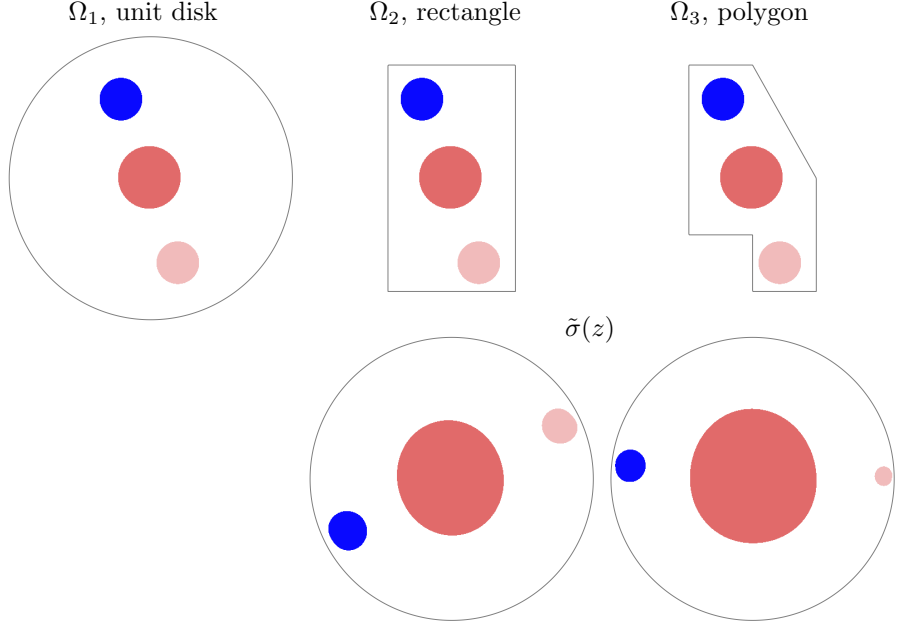

\subsection{Reconstructions without magnification}
\label{sec:eka}
We consider three different {\em true domains}:  
\begin{enumerate}
\item $\Omega_1 = D$ is the unit disk,
\item $\Omega_2\subset D$ is a rectangle, 
\item $\Omega_3\subset D$ is a more complicated polygonal domain.
\end{enumerate}
The target conductivity $\sigma_1 = \sigma: D \to \R_+$ for $\Omega_1$ is composed of three embedded inclusions in homogeneous unit background; see the top left image of Figure~\ref{fig:domains}. The red inclusion in the middle has constant conductivity 3, whereas the blue and light red ones closer to the boundary are characterized by conductivity levels 0.2 and 2, respectively. The target conductivities for $\Omega_2$ and $\Omega_3$ are defined as $\sigma_2 = \sigma|_{\Omega_2}$ and $\sigma_3 = \sigma|_{\Omega_3}$, respectively; observe that ${\rm supp} (\sigma - 1) \subset \Omega_j$, $j= 2,3$, so that the target conductivity equals one in some interior neighborhood of the boundary for all three target domains.
The bottom row of Figure~\ref{fig:domains} illustrates the conformally mapped conductivities $\tilde{\sigma}_j = \tilde{\sigma}_j \circ \Psi_j$, $j=2,3$, in the unit disk. Here, $\Psi_j$, $j=2, 3$, is a Schwarz--Christoffel map sending the unit disk onto $\Omega_j$, $j=2,3$, respectively, with the extra requirement $\Psi_j(0) = 0$ that fixes the map up to rotations of the unit disk.

\begin{figure}
\begin{picture}(120,85)
\epsfxsize=4cm
\put(0,40){\includegraphics*[width=3.8cm]{\imagepath 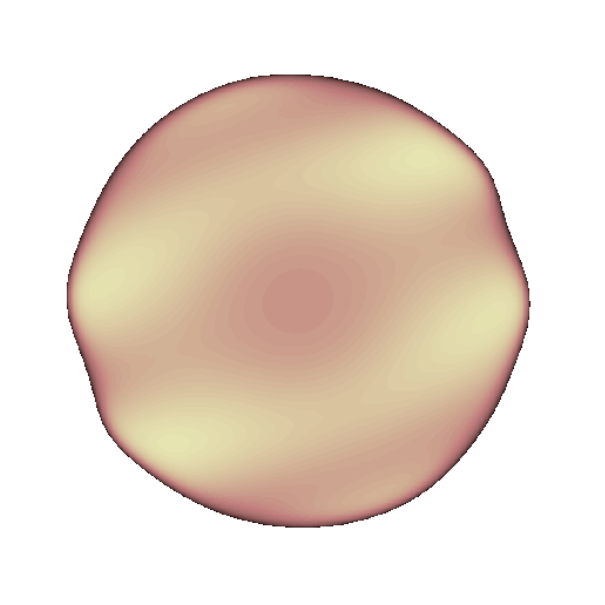}}
\put(40,40){\includegraphics*[width=3.8cm]{\imagepath 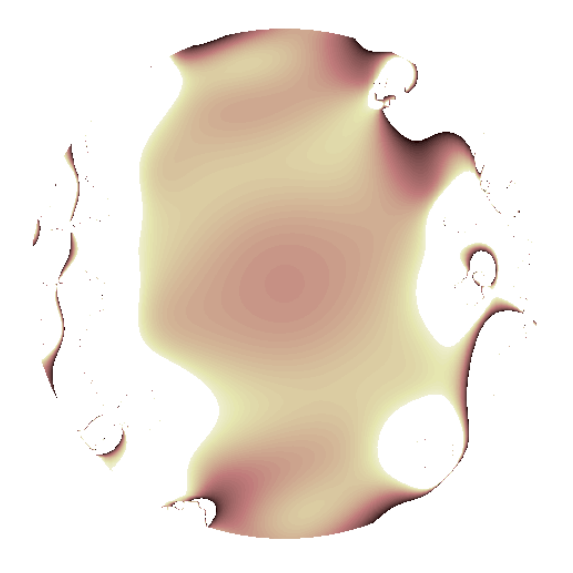}}
\put(80,40){\includegraphics*[width=3.8cm]{\imagepath 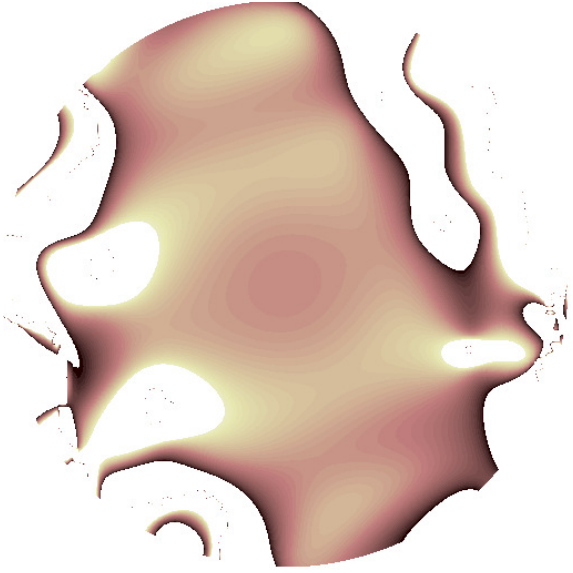}}

\put(40,0){\includegraphics*[width=3.8cm]{\imagepath 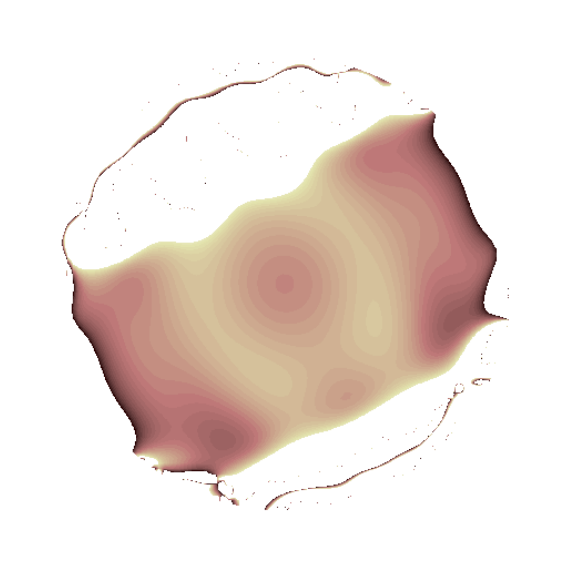}}
\put(80,0){\includegraphics*[width=3.8cm]{\imagepath 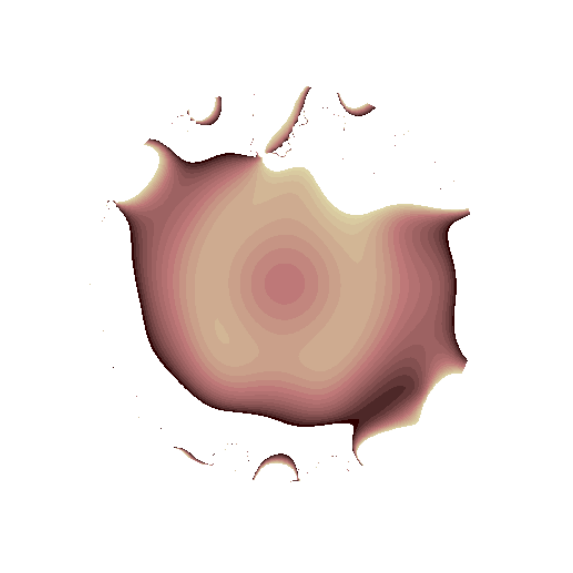}}
\put(8,80){\small $\Omega_1$, unit disk}
\put(48,80){\small $\Omega_2$, rectangle}
\put(88,80){\small $\Omega_3$, polygon}
\put(60,35){\small Conformally transformed}
\put(28,40){$R = 8$}
\put(65,40){$R = 9$}
\put(105,40){$R = 10$}
\put(65,0){$R = 8$}
\put(105,0){$R = 7$}
\end{picture}
\caption{Numerical approximations for the real part of the truncated scattering transform $\T_{R,c}$ of \eqref{scatRc} for the five configurations in Figure~\ref{fig:domains} with $c = 10$. 
The white parts correspond to $\abs{\T}>c$ or $\abs{k}>R$,~i.e,~$\T_{R,c}=0$. 
\label{fig:scats}}
\end{figure}
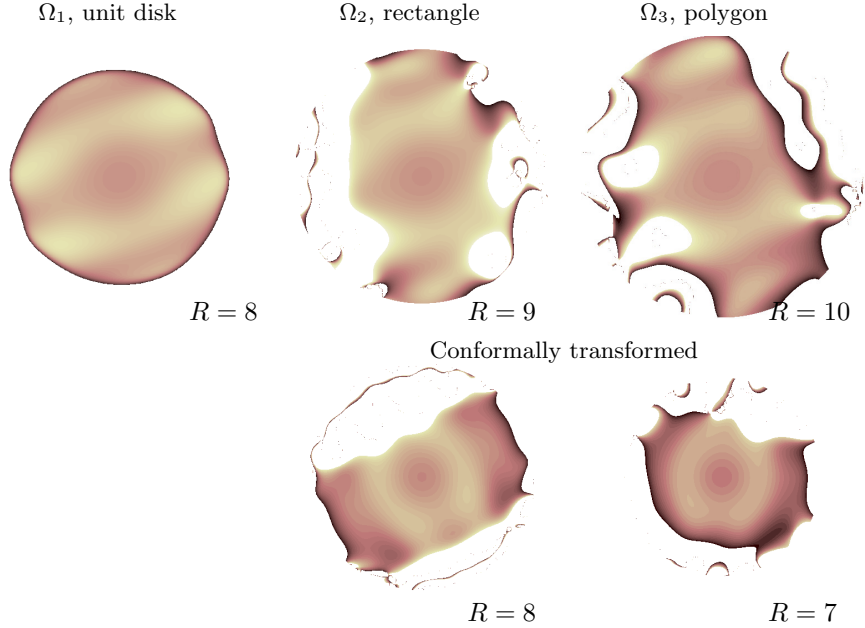

\begin{figure}
\begin{picture}(120,125)
\put(0,80){\includegraphics*[width=3.8cm]{\imagepath EITconfV7_disk_sigma-eps-converted-to.pdf}}
\put(40,80){\includegraphics*[width=3.8cm]{\imagepath EITconfV7_rect_sigma-eps-converted-to.pdf}}
\put(80,80){\includegraphics*[width=3.8cm]{\imagepath EITconfV7_polygon_sigma-eps-converted-to.pdf}}

\put(0, 40){\includegraphics*[width=3.8cm]{\imagepath 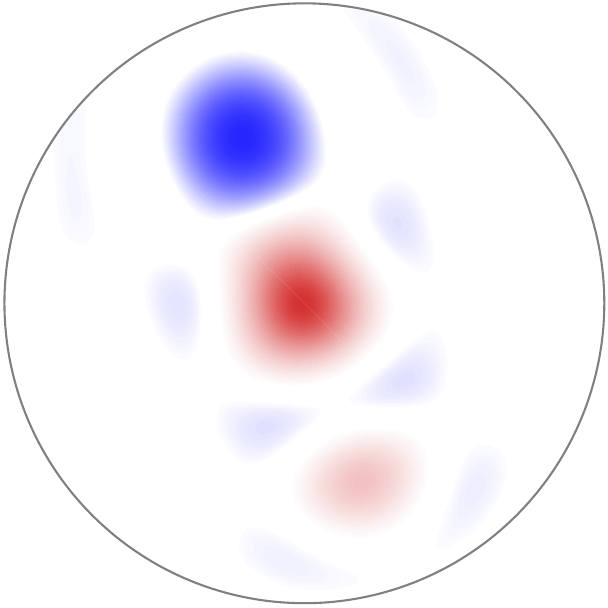}}
\put(40,40){\includegraphics*[width=3.8cm]{\imagepath 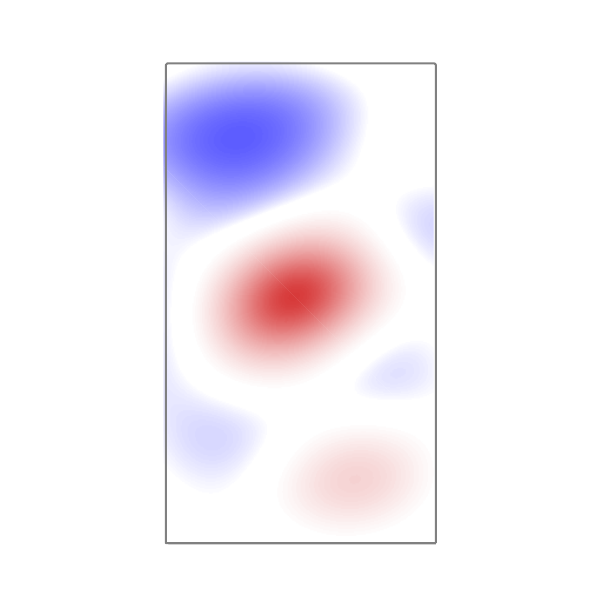}}
\put(80,40){\includegraphics*[width=3.8cm]{\imagepath 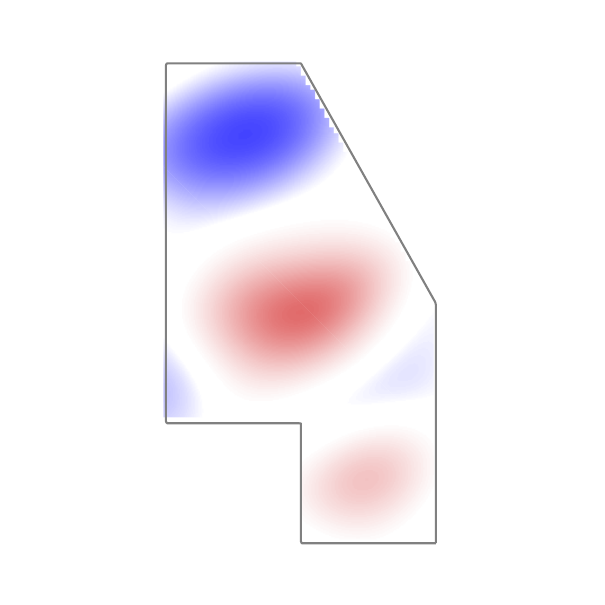}}

\put(40,0){\includegraphics*[width=3.8cm]{\imagepath 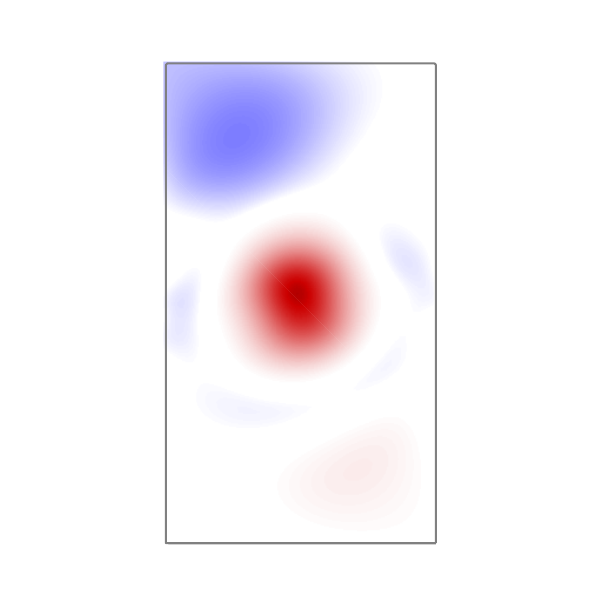}}
\put(80,0){\includegraphics*[width=3.8cm]{\imagepath 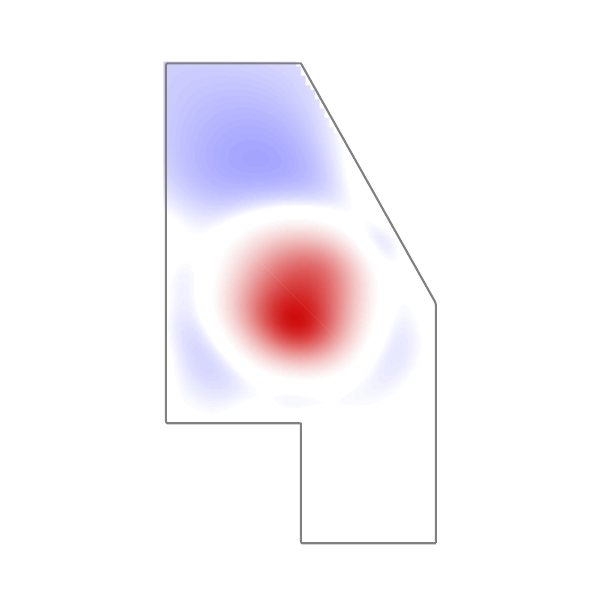}}

\put(8,120){\small Target in $\Omega_1$}
\put(48,120){\small Target in $\Omega_2$}
\put(88,120){\small Target in $\Omega_3$}
\put(52,77){\small Standard D-bar}
\put(52,37){\small Conformal D-bar}

\put(30,42){\small 18.9\%}
\put(67,42){\small 25.7\%}
\put(103,42){\small 29.2\%}
\put(67,3){\small 26.2\%}
\put(103,3){\small 27.9\%}

\end{picture}
\caption{Reconstructions of the test conductivity in three different domains and  the corresponding relative $L^2(\Omega)$-errors. Top row: the target configurations. Middle row: reconstructions by the D-bar method using the current patterns \eqref{Fbasis} with $N=16$. Bottom row: reconstructions by the D-bar method in the virtual domain (i.e., the unit disk) using the conformally transformed current patterns \eqref{eq:trans} with $N=16$ in the true domain.
\label{fig:recons}}
\end{figure}

Figure~\ref{fig:scats} shows the real parts of the truncated scattering transforms \eqref{scatDN2} for the five geometries of Figure~\ref{fig:domains} with $c=10$.
The reconstructions produced by the D-bar method are illustrated in Figure~\ref{fig:recons}, with the top row showing the target conductivities for comparison. The reconstructions formed in the original domains using the respective standard Fourier current basis \eqref{Fbasis} are shown in the middle row. The bottom row presents the reconstructions corresponding to the employment of the conformally transformed Fourier current basis \eqref{eq:trans} on $\partial \Omega_j$, $j=2,3$, subsequently applying the unit disk version of the D-bar method (see~Section~\ref{sec:recon}) to the resulting boundary data (cf.~Lemma~\eqref{changevar}) to reconstruct the virtual conductivity $\tilde{\sigma}_j$, $j=2,3$, in $D$, and finally mapping the reconstruction back to the true domain with the help of $\Psi_j$, $j=2,3$. The reconstructions in the second and third rows emphasize different characteristic of the target conductivities. The (single) reconstruction in $\Omega_1$ is arguably the best one, but for $\Omega_2$ and $\Omega_3$ the two reconstruction techniques seem to function approximately as well. Hence, it seems that the choice between the two of them should be made by comparing the respective computational cost, that is, whether the possibility of implementing the D-bar method in the unit disk (see~Section~\ref{sec:recon}) compensates for the extra cost of forming the conformal map between the true domain and the unit disk. The answer to this question depends on many things such as,~e.g.,~whether the reference N-to-D map corresponding to unit conductivity for the true domain can measured or needs to be computed and whether there exists an efficient technique for evaluating the needed conformal map for the considered measurement geometry.

\subsection{Magnification of a ROI}
\label{sec:toka}
\begin{figure}
\begin{picture}(120,130)
\put(100,62){\includegraphics*[width=1cm]{\imagepath 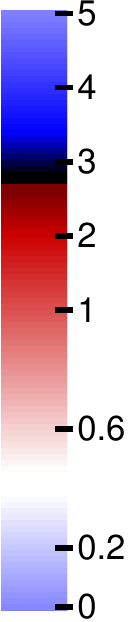}}
\put(0,86){\includegraphics*[width=4cm]{\imagepath 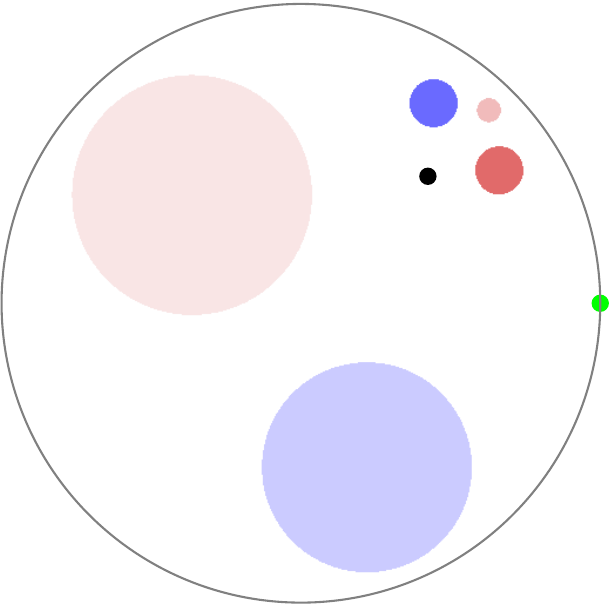}}
\put(50,86){\includegraphics*[width=4cm]{\imagepath 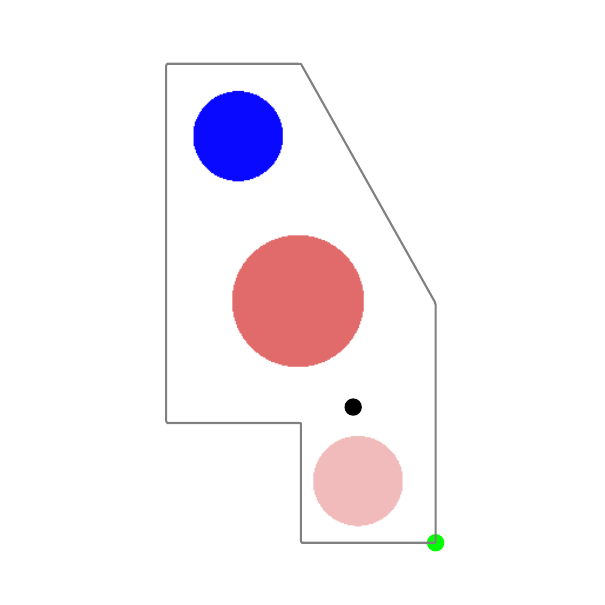}}
\put(0,44){\includegraphics*[width=4cm]{\imagepath 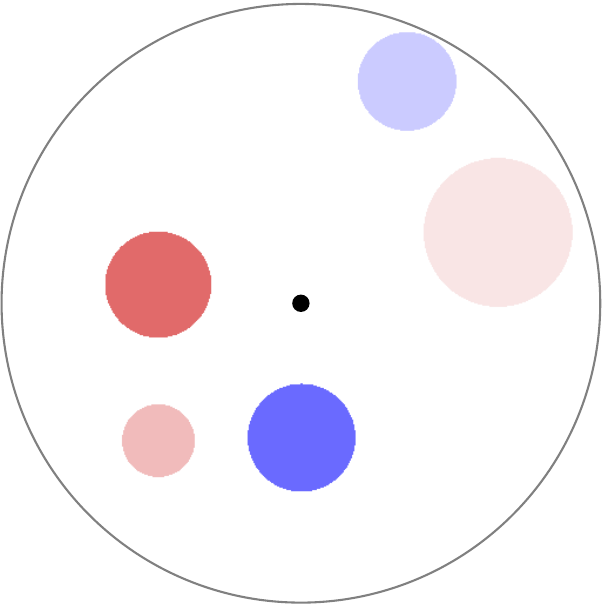}}
\put(50,44){\includegraphics*[width=4cm]{\imagepath 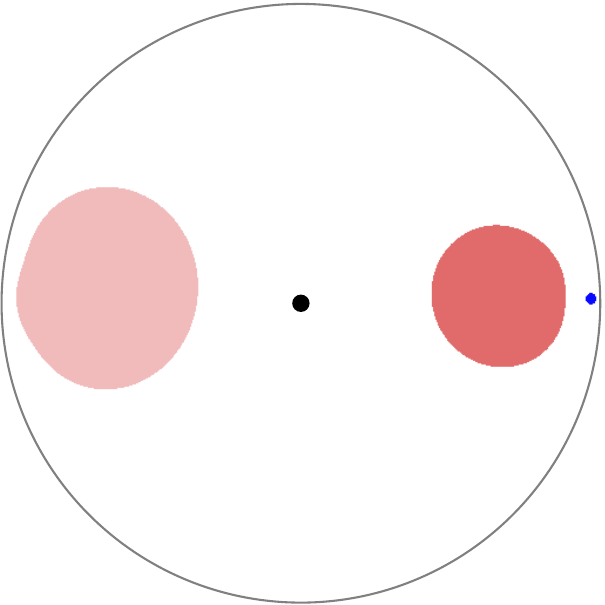}}
\put(0,0){\includegraphics*[width=4cm]{\imagepath 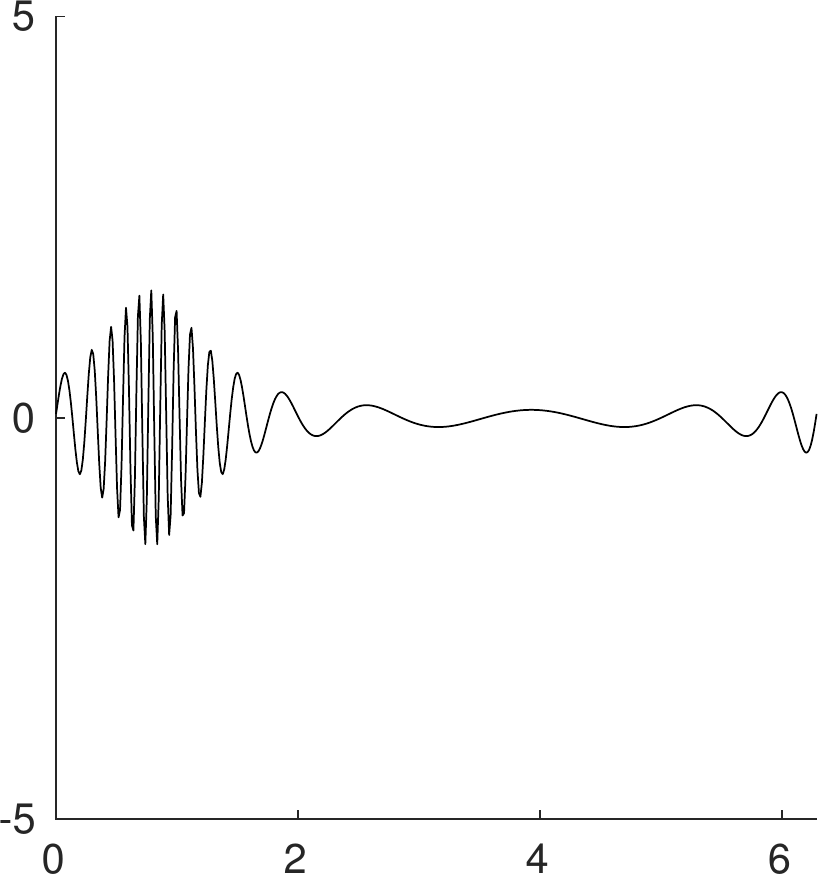}}
\put(50,0){\includegraphics*[width=4cm]{\imagepath 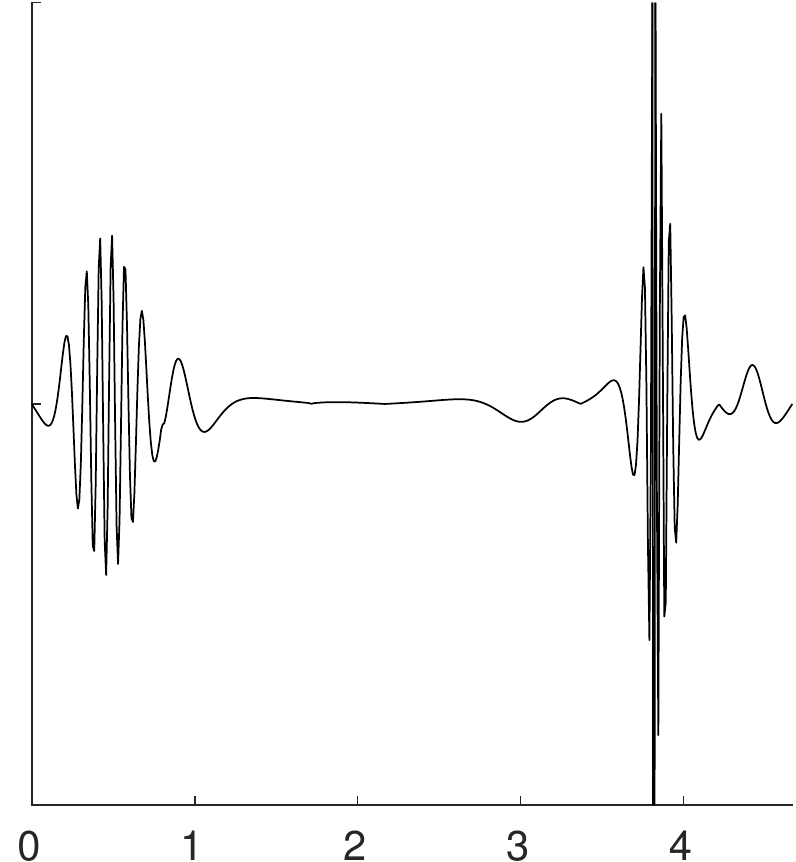}}
\end{picture}
\caption{Magnification of a ROI. Top row: the true target configurations. 
Middle row: the conformally transformed conductivities in the virtual domain with magnified ROIs. On their boundaries the point corresponding to $s=0$ is marked with a green dot. The black dot is in the ROI and mapped to the origin. Bottom row: The real parts of the conformally transformed Fourier current pattern $\varphi_{16}$ of \eqref{eq:trans} for the respective true target domains.
\label{fig:MobiusDomains}}
\end{figure}

In this section, we test the idea of magnifying a ROI by a conformal map. The two target conductivities and true domains are shown in the top row of Figure~\ref{fig:MobiusDomains}. The left-hand domain is the unit disk with several embedded circular inhomogeneities; the first quadrant with the three small inclusions is considered as the ROI. The right-hand domain is the polygon $\Omega_3$ from the previous section with the same conductivity pattern consisting of three inclusions in a homogeneous background; the bottom antler enclosing the pink inclusion serves as the ROI. 

The middle row of Figure~\ref{fig:MobiusDomains} shows two virtual conductivities in the unit disk obtained by magnifying the ROI with a certain conformal map; in both cases, the requirement that the black dot in the top row is mapped to the corresponding one in the middle row fixes $\Phi$ (up to rotations of the unit disk). For the left-hand setup, such $\Phi$ is realized by the M\"obius transform $\MM_a$ of \eqref{mobius} with the free parameter $a = 0.6\exp(\rmi \theta/4)$. For the right-hand configuration $\Phi$ is a Schwarz--Christoffel map, which we construct by composing  a suitable M\"obius transform with a preliminary Schwarz--Christoffel map that sends the origin to itself. The highest,~i.e.~the sixteenth, conformally transformed Fourier current pattern (cf.~\eqref{eq:trans}) for each domain is shown in the bottom row of Figure~\ref{fig:MobiusDomains}. Driving these current densities through the boundaries of the true domains corresponds to the use of the standard Fourier boundary current $\tilde{\phi}_{16}$ for the virtual configurations in the middle row of Figure~\ref{fig:MobiusDomains}.

\begin{figure}
\begin{picture}(120,135)
\put(110,40){\includegraphics*[width=1cm]{\imagepath EITconfV7_colorbar_arXiv.pdf}}
\put(0,89){\includegraphics*[width=4cm]{\imagepath 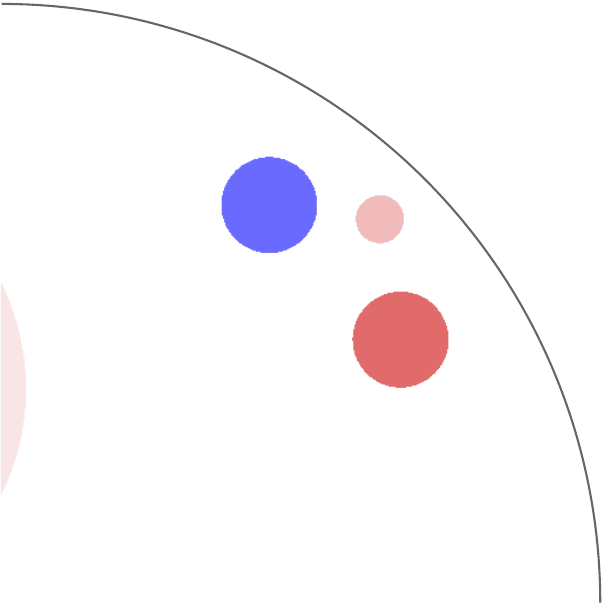}}
\put(57,89){\includegraphics*[height=4cm]{\imagepath 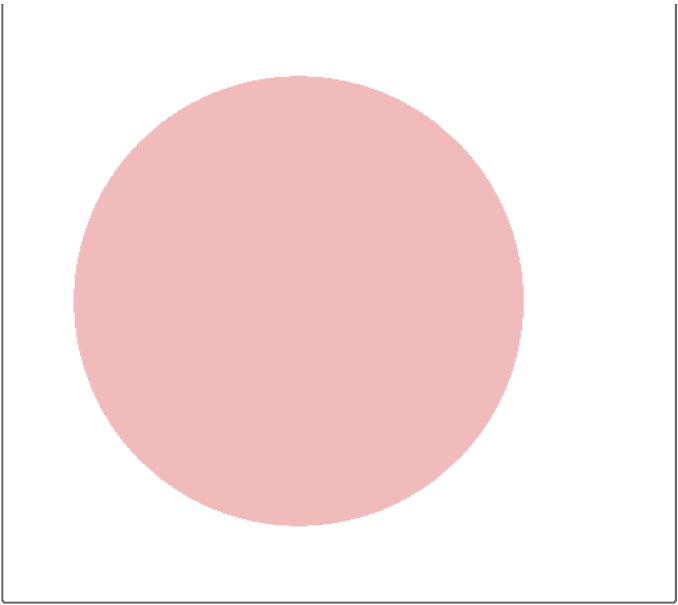}}
\put(0,45){\includegraphics*[width=4cm]{\imagepath 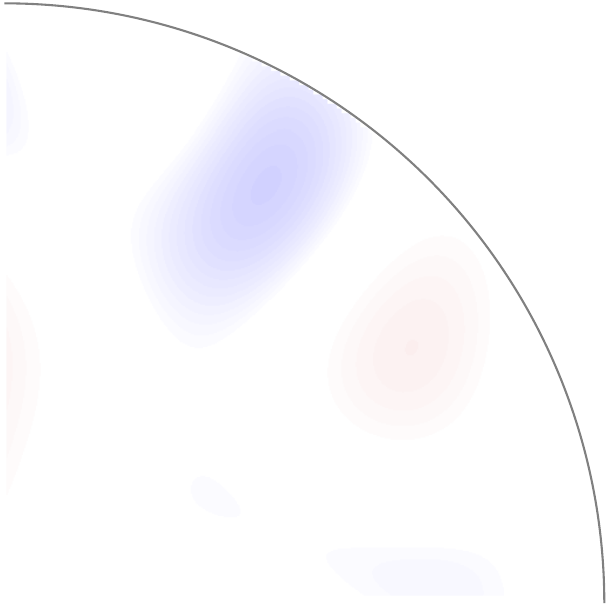}}
\put(57,45){\includegraphics*[height=4cm]{\imagepath 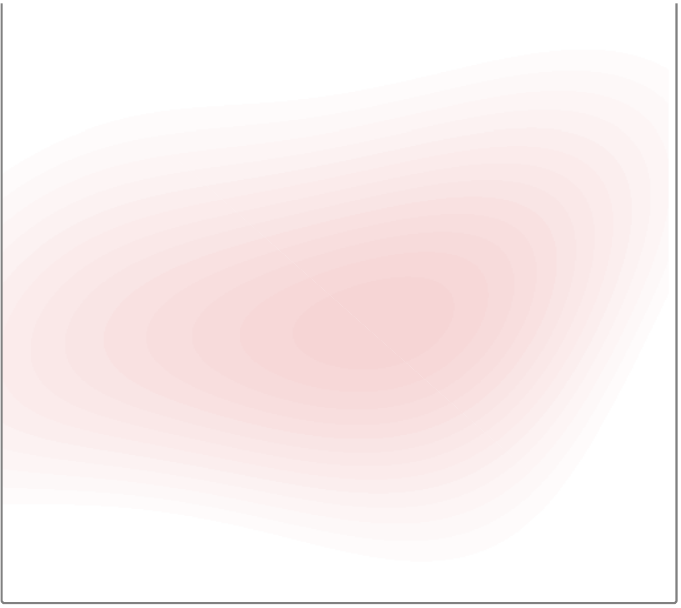}}
\put(0,0){\includegraphics*[width=4cm]{\imagepath 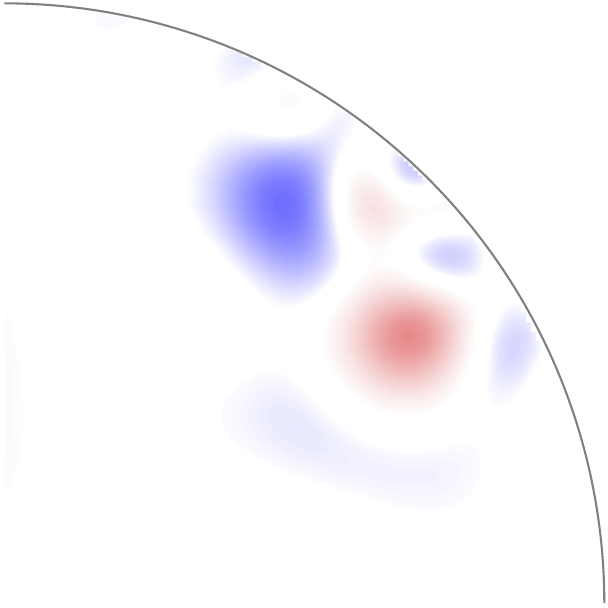}}
\put(57,0){\includegraphics*[height=4cm]{\imagepath 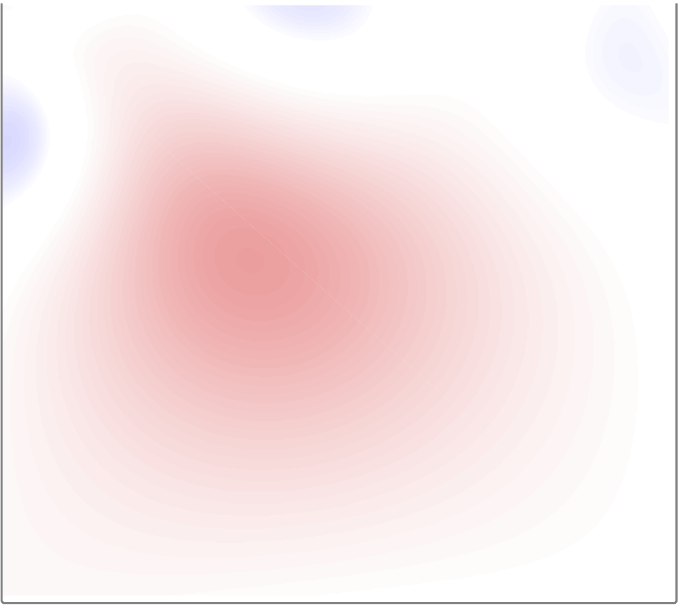}}

\put(22,125){\small Target conductivity}
\put(33,121){\small  in the ROI}

\put(28,80){\small Standard D-bar}
\put(29,35){\small ROI-magnifying D-bar}

\put(34,46){\small 27.3\%}
\put(34,1){\small 18.1\%}
\put(90,46){\small 26.5\%}
\put(90,1){\small 20.5\%}

\end{picture}
\caption{Reconstructions by magnifying the ROIs (cf.~Figure~\ref{fig:MobiusDomains}) with the corresponding relative $L^2$-errors over the ROIs. Top row: the ROIs in the original domains of Figure~\ref{fig:MobiusDomains}. 
Middle row: reconstructions by the D-bar method using the current patterns \eqref{Fbasis} with $N=16$. Bottom row: reconstructions by the D-bar method in the virtual domain (i.e., the unit disk) using the conformally transformed current patterns \eqref{eq:trans} with $N=16$ in the true domain.
\label{fig:MobiusRecons}}
\end{figure}

The first row of Figure~\ref{fig:MobiusRecons} shows the target conductivities within the ROIs. The conductivity reconstructions produced by the D-bar method implemented directly in the true target domains with the standard Fourier current patterns \eqref{Fbasis} are presented in the middle row. The bottom row illustrates the reconstructions corresponding to the use of conformally transformed Fourier currents \eqref{eq:trans} for the true domain and an application of the standard D-bar method in the virtual reconstruction domain (i.e.,~the unit disk). The reconstructions in Figure~\ref{fig:MobiusRecons} demonstrate that the magnification by a suitable conformal map 
results in a more accurate reconstruction of the conductivity within the ROI. For all reconstructions in Figure~\ref{fig:MobiusRecons}, we set $c=20$ in \eqref{scatDN2}, while $R>0$ was chosen as explained just before Section~\ref{sec:eka}.

\subsection{Reconstructions in the half-plane}
\label{sec:halfspacenumerics}
Our final numerical experiment considers the unbounded case $\Omega = \{ z \in \C \ | \ {\rm Im}(z) > 0 \}$; see Appendix~\ref{sec:halfplane} for information on the well posedness of the forward problem of EIT in such a half-space geometry as well as on the employed layer potential techniques for simulating the needed boundary data. The M\"obius transforms used for mapping $\Omega$ onto the unit disk are of the form
$$
\MM_{b}(z) = \frac{(z-b)-\xi}{(z-b)-\kk{\xi}},
$$
where $\xi = 1.2\rmi$ is fixed and the horizontal translation parameter $b \in \R$ takes five different values, namely $b = -3, -1.5, 0, 1.5, 3$. In particular, $\MM_b$ maps the point $b+\xi$ to the origin.

\begin{figure}
\begin{picture}(120,125)
\put(0,105){\includegraphics*[width=8cm]{\imagepath 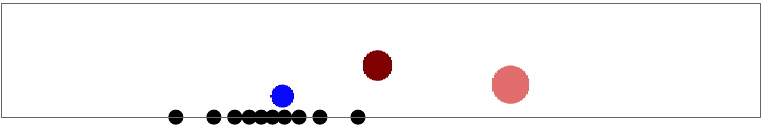}}
\put(0,80){\includegraphics*[width=8cm]{\imagepath 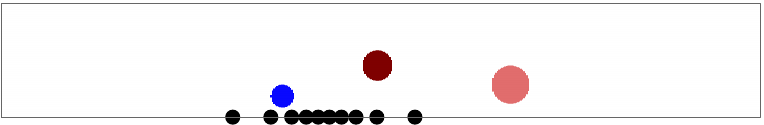}}
\put(0,55){\includegraphics*[width=8cm]{\imagepath 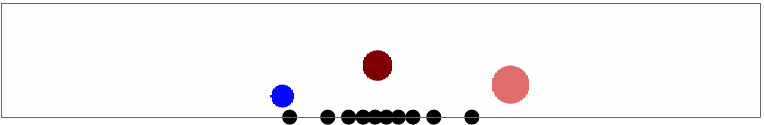}}
\put(0,30){\includegraphics*[width=8cm]{\imagepath 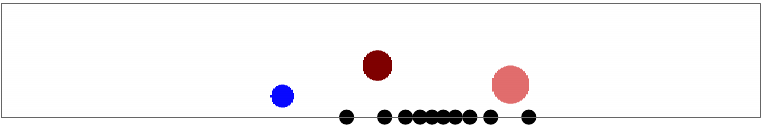}}
\put(0,5){\includegraphics*[width=8cm]{\imagepath 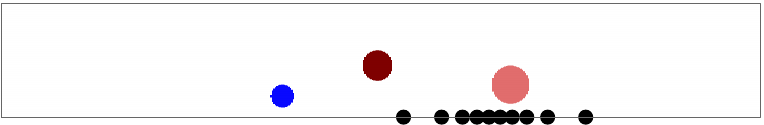}}
\put(85,100){\includegraphics*[height=2.5cm]{\imagepath 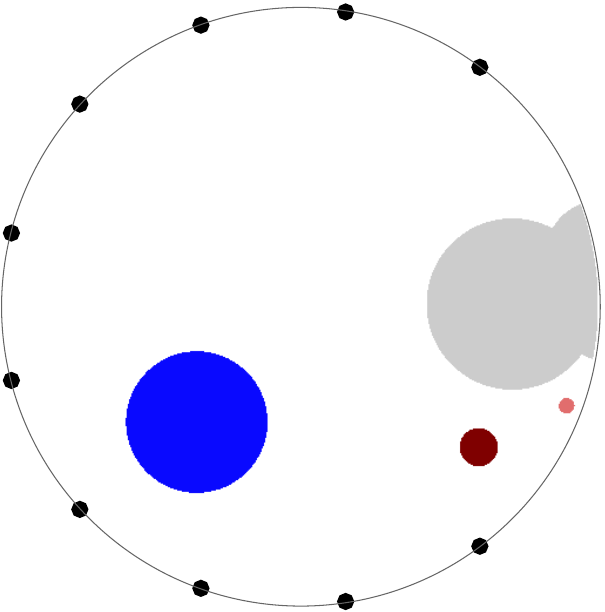}}
\put(85,75){\includegraphics*[height=2.5cm]{\imagepath 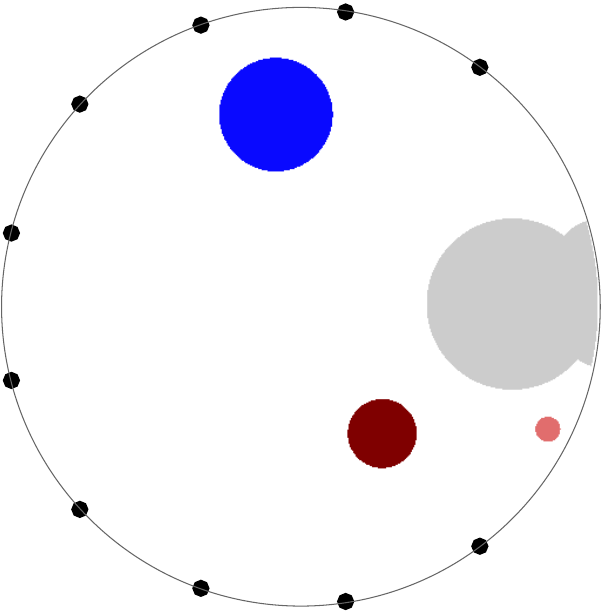}}
\put(85,50){\includegraphics*[height=2.5cm]{\imagepath 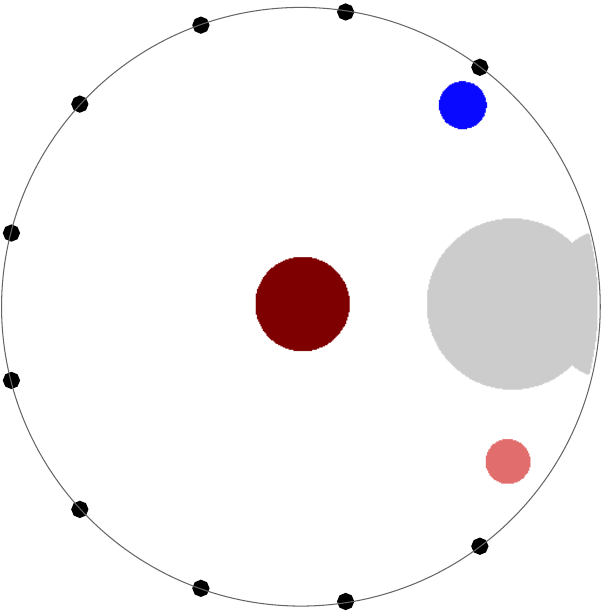}}
\put(85,25){\includegraphics*[height=2.5cm]{\imagepath 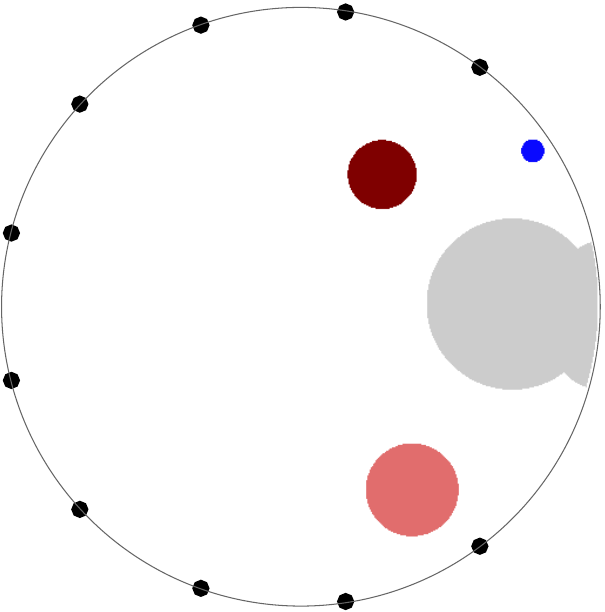}}
\put(85,0){\includegraphics*[height=2.5cm]{\imagepath 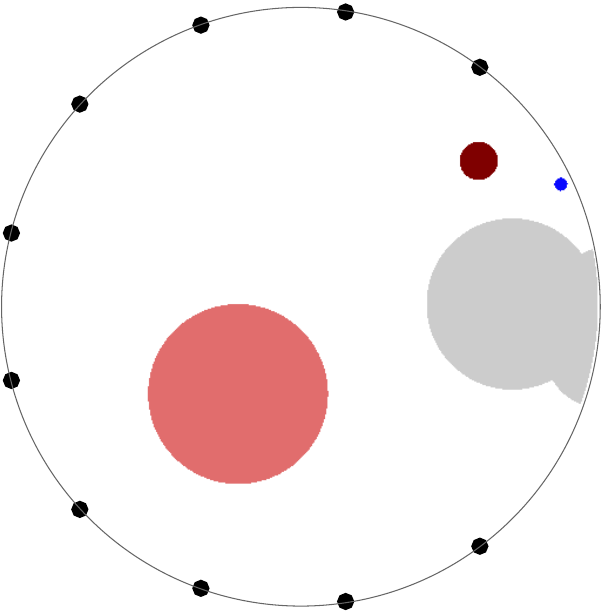}}
\put(30,101){$b = -3$}
\put(30,76){$b = -1.5$}
\put(30,51){$b = 0$}
\put(30,26){$b = 1.5$}
\put(30,1){$b = 3$}
\end{picture}
\caption{Conductivity inhomogeneities in the upper half-plane (left) and the corresponding virtual conductivities in the unit disk corresponding to $\MM_b$ with five different values for $b$ (right). The black dots illustrate an electrode array that could be used for approximating the conformally transformed Fourier currents \eqref{eq:trans} with $\Phi = \MM_b$ and $N=5$. The gray areas on the right correspond to the exterior of the rectangular subset considered on the left. \label{fig:halfplanedemo}}
\end{figure}

The target conductivity consists of three inclusions in homogeneous unit background; it is visualized in the left-hand column of Figure~\ref{fig:halfplanedemo} over the subset $\Omega_0 = \{z\in\C \ | \ \abs{\re(z)}\leq 10, \ 0\leq\im(z)\leq 3\} \subset \Omega$. The right-hand column of Figure~\ref{fig:halfplanedemo} shows the virtual conductivities $\tilde{\sigma} = \sigma \circ \Psi$ in the unit disk for the conformal map $\Phi = \MM_{b}$ and our five choices of $b \in \R$. The gray area represents $\MM_{b}(\Omega \setminus \Omega_0)$. Observe that the parameter $b$ determines the horizontal position that is magnified the most under $\MM_{b}$ (or, actually, shrunk the least). 

\begin{figure}
\begin{picture}(120,120)
\put(114,35){\includegraphics*[width=1cm]{\imagepath EITconfV7_colorbar_arXiv.pdf}}
\put(-1,100){\includegraphics*[width=11cm]{\imagepath 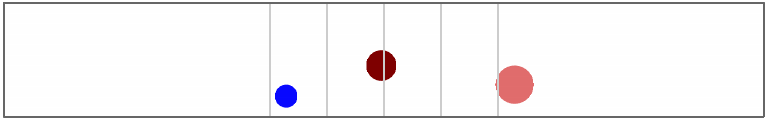}}
\put(-1,80){\includegraphics*[width=11cm]{\imagepath 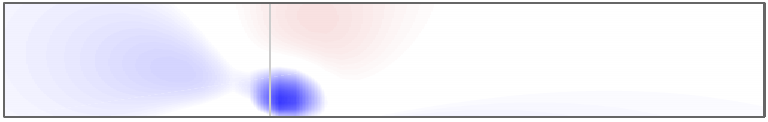}}
\put(-1,60){\includegraphics*[width=11cm]{\imagepath 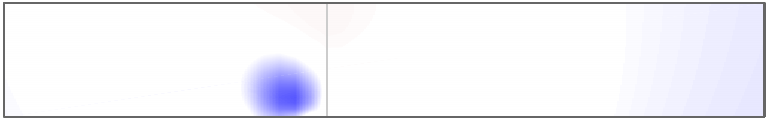}}
\put(-1,40){\includegraphics*[width=11cm]{\imagepath 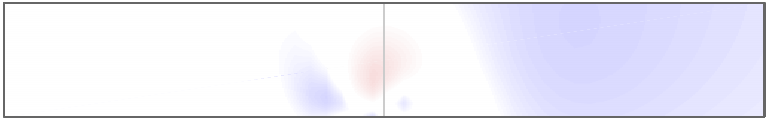}}
\put(-1,20){\includegraphics*[width=11cm]{\imagepath 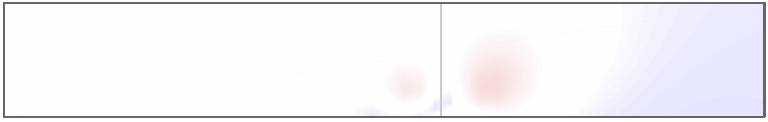}}
\put(-1,0){\includegraphics*[width=11cm]{\imagepath 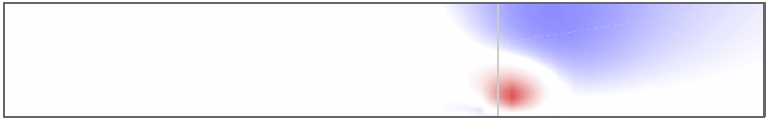}}
\put(2,82){$b = -3$}
\put(2,62){$b = -1.5$}
\put(2,42){$b = 0$}
\put(2,22){$b = 1.5$}
\put(2,2){$b = 3$}
\end{picture}
\caption{D-bar reconstruction in the upper half-plane obtained by using five different M\"obius transforms,~i.e.,~$\MM_b$ with $b=-3,-1.5,0,1.5,3$, for mapping the half-plane onto the unit disk. The vertical line segments indicate the used values of $b$ (cf.~Figure~\ref{fig:halfplanedemo}).\label{fig:halfplanerecon}}
\end{figure}

For this half-plane setting, we use only ten conformally transformed Fourier currents on $\partial \Omega$,~i.e.,~set $N=5$ in \eqref{eq:trans}.
The black dots on the boundary of $\Omega$ in Figure~\ref{fig:halfplanedemo} represent a possible configuration of small electrodes for approximating the conformally transformed Fourier currents~\eqref{eq:trans} corresponding to the particular $\Phi = \MM_{b}$; the images of these `electrodes' under $\MM_{b}$ form an (incomplete) uniform grid on the boundary of the virtual domain as shown in the right-hand column of Figure~\ref{fig:halfplanedemo} (cf.~Appendix~\ref{sec:point} and \ref{sec:halfplane}). To be more precise, for each $\Phi = \MM_{b}$ the parts of (the PEM approximation for) the current patterns \eqref{eq:trans} supported outside the interval $[-2.5+b,2.5+b]$ on the real axis are ignored,~i.e.,~set to zero, and the corresponding potential measurements are also only taken on this interval and set to zero on the rest of $\partial \Omega$. The effect of these truncations is not further analyzed, but increasing the width of the measurement interval would certainly improve the reconstructions to a certain extent --- as would increasing the number of employed conformally transformed Fourier currents on $\partial \Omega$.

The reconstructions are once again formed by the D-bar method in the virtual domain (i.e., the unit disk) and mapped back to the half-plane by $\MM_{b}^{-1}$.  
The five reconstructions corresponding to our five values for the parameter $b$ are shown in Figure~\ref{fig:halfplanerecon}, with the vertical gray lines indicating the midpoints of the `measurement arrays',~i.e.,~the five values of $b$. It is obvious that the parameter $b$ pinpoints the area that is reconstructed most accurately: The blue and light red inclusions are clearly visible in the first ($b=-3$) and the last ($b=3$) reconstruction, respectively. In addition, there are arguably some traces of the dark red inclusion, lying the furthest from $\partial \Omega$, in the third reconstruction ($b=0$). 
One could postprocess the five reconstructions of Figure~\ref{fig:halfplanerecon} to obtain a `composite reconstruction', but such image processing is outside the scope of this work.

\section{Concluding remarks}\label{sec:conclusions}
By utilizing conformal maps, we have introduced and numerically tested a method for transferring the D-bar method from an arbitrary simply connected planar domain to the unit disk, where the algorithm can be implemented most efficiently. In particular, we have presented the first D-bar reconstructions in a half-plane. By exploiting the degrees of freedom in the choice of a conformal map between a given simply connected planar domain and the unit disk (characterized by the M\"obius transforms sending the unit disk onto itself), we have also introduced a method for magnifying a ROI and thus obtaining more accurate reconstructions of its conductivity. The magnification of a ROI is based on the use of so-called conformally transformed Fourier current patterns that exhibit fast spatial oscillations close to the ROI. However, we have proved that idealized continuum boundary data  corresponding to such boundary current densities can be accurately approximated by real-world EIT measurements if the available (small) electrodes are positioned in a suitable nonuniform configuration; see Appendix~\ref{sec:point}.

\appendix

\section{Approximations by point electrodes}
\label{sec:point}

The aim of this appendix is to explain how the positioning of the available electrodes on the boundary of the true domain is related to the ability to approximate the output of the relative N-to-D map for given virtual boundary current densities on the boundary of the unit disk. To put it short, the electrode pattern should be densest on those boundary sections where the  conformally transformed current densities oscillate the most. 

For simplicity, we assume in the following that $\partial \Omega$ is of the class $C^\infty$, although the considered distributional Neumann boundary value problems could be defined even in domains with piecewise $C^{1,\alpha}$ boundaries~\cite{Seiskari14}. As in Section~\ref{sec:conformal}, $\Phi: \Omega \to D$, with $\Psi = \Phi^{-1}$, is a conformal map of the true domain $\Omega$ onto the virtual domain,~i.e.,~the open unit disk. In addition, it is assumed that the conductivity $\sigma$ equals one in some interior neighborhood of $\partial \Omega$. Under these assumptions, the relative N-to-D map
\begin{equation}
\label{eq:boundedness}
\RR_\sigma - \RR_1: H_\diamond^r(\partial \Omega) \to H^{q}(\partial \Omega)/\C 
\end{equation}
is bounded for any $r,q \in \R$~\cite{Lions72}. As $\tilde{\sigma} = \sigma\circ \Psi$ also equals one in some interior neighborhood of $\partial D$, the analogous result obviously holds for $\tilde{\RR}_{\tilde{\sigma}} - \tilde{\RR}_1$ as well, with $\tilde{\RR}_{\tilde{\sigma}}, \tilde{\RR}_1$ being N-to-D maps for the unit disk.

The PEM models the electrodes as pointlike boundary current sources. Let $\tilde{z}_m \in \partial D$, $m=0, \dots, M-1$, be a set of {\em virtual} equiangular electrodes on the boundary of the unit disk and define the corresponding {\em true} point electrodes on $\partial \Omega$ to be $z_m = \Psi(\tilde{z}_m)$, $m=0, \dots, M-1$. Hence, the point electrodes for the true domain are the preimages of equidistant electrodes on $\partial D$ under the bijective and smooth map $\Phi: \overline{\Omega} \to \overline{D}$. These fixed electrode positions are used in the following notation: for arbitrary $g\in H^r(\partial \Omega)$ and $\tilde{g}\in H^r(\partial D)$, $r>1/2$, we write
\begin{align*}
[g]_M = [g(z_m)]_{m=0}^{M-1} \in\C^M \qquad {\rm and} \qquad
[\tilde{g}]_M = [\tilde{g}(\tilde{z}_m)]_{m=0}^{M-1} \in\C^M.
\end{align*}
In particular, if $g = \tilde{g}\circ\Phi|_{\partial \Omega}$, then $[g]_M = [\tilde{g}]_M$. 

Let $\C_\diamond^M$ be the zero-mean subspace of $\C^M$ and denote by $\I \in \C_\diamond^M$ a net current pattern through the pointlike electrodes $\{z_m\}_{m=0}^{M-1}$ on the true domain boundary $\partial \Omega$.  The relative current-to-potential operator $\MMM: \C_\diamond^M \to \C^M/\C$ of the PEM is defined via~\cite{Hanke2011}
\begin{equation}
\label{eq:Adef}
\MMM: \I  \mapsto \left[\left((\RR_\sigma - \RR_1)\Big(\sum_{m=0}^{M-1} \I_m \delta_{z_m}\Big) \right) \!\right]_M,
\end{equation}
where $\delta_{z} \in H^{-1/2 -\epsilon}(\partial \Omega)$, $\epsilon >0$, is the Dirac delta distribution supported at $z \in \partial \Omega$. In other words, the net electrode currents $\I$ are driven through the pointlike electrodes and the corresponding relative potentials are measured at the same locations. Notice that $\MMM$ is well-defined due to \eqref{eq:boundedness} and the Sobolev embedding theorem. From the practical standpoint, it is important to note that the discrepancy between $\MMM$ and the relative current-to-potential map of the CEM \cite{Cheng1989,Somersalo1992} is of the order $O(d^2)$, with $d>0$ being the maximal diameter of the electrodes~\cite{Hanke2011}.

Suppose a given {\em virtual} continuum current pattern $\tilde{f} \in H^{r}_\diamond(\partial D)$, $r > 1/2$, on the boundary of the unit disk satisfies $[\tilde{f}]_M\in \C_\diamond^M$. In particular, this holds if $\tilde{f}$ is an element of the Fourier current basis \eqref{Fbasis}, with $N<M$, on $\partial D$ since $\tilde{z}_m$, $m=0, \dots, M-1$, are equiangular. We define the corresponding {\em true} electrode net currents as
\begin{equation}
\label{eq:Ifcurrent}
\I^{\tilde{f}} := \frac{2\pi}{M} [\tilde{f}\circ\Phi]_M = \frac{2\pi}{M} [\tilde{f}]_M \in \C_\diamond^M.
\end{equation}
The following theorem demonstrates that the virtual relative potential measurement on $\partial D$ corresponding to $\tilde{f}$ can be approximated by the true electrode measurements that correspond to driving the net currents $\I^{\tilde{f}} \in \C_\diamond^M$ through the point electrodes $\{z_m\}_{m=0}^{M-1} \subset \partial \Omega$.

\begin{Thm} 
\label{main_theorem}
Assume that $\tilde{f} \in H^r_\diamond(\partial D)$, $r\in \N$ and $[\tilde{f}]_M\in \C_\diamond^M$. Then, it holds that 
\begin{equation}
\label{eq:optconv}
\Big\| \MMM(\I^{\tilde{f}}) - \big[\big((\tilde{\RR}_{\tilde{\sigma}} - \tilde{\RR}_1)\tilde{f}\big)\big]_M \Big\|_{\C^M/\C} \leq \frac{C}{M^{r-1/2}} \| \tilde{f} \|_{H^r(\partial D)},
\end{equation}
where $C = C(r,\tilde{\sigma})>0$ is independent of $M\in \N \setminus \{1\}$ and $\tilde{f}$.
\end{Thm}

\begin{proof}
First of all, it follows from the compatibility of conformal maps and Neumann boundary values composed of Dirac deltas that
\begin{equation}
\label{eq:deltacompo}
\left((\RR_\sigma - \RR_1)\Big(\sum_{m=0}^{M-1} \I_m^{\tilde{f}} \delta_{z_m}\Big)\right)\circ \Psi|_{\partial D}
= (\tilde{\RR}_{\tilde{\sigma}} - \tilde{\RR}_1)\Big(\sum_{m=0}^{M-1} \I_m^{\tilde{f}} \delta_{\tilde{z}_m}\Big),
\end{equation}
which is to be understood modulo an additive constant; see the proof of \cite[Theorem~3.2]{Hakula11}. To shorten the notation, we denote
$$
\tilde{f}_M = \sum_{m=0}^{M-1} \I_m^{\tilde{f}} \delta_{\tilde{z}_m} \in H^{-q}_\diamond(\partial D), \qquad q > 1/2,
$$
in what follows.

Let $r \in \N$ and $v \in H^r(\partial D)/\C$ be arbitrary, and denote by $v_0 \in H^r_\diamond(\partial D)$ the zero-mean element in the equivalence class $v \in H^r(\partial D)/\C$. In particular, take note that $\|v \|_{H^r(\partial D)/\C} = \|v_0 \|_{H^r(\partial D)}$. Since the mean-free distribution $\tilde{f}_M - \tilde{f} \in H^{-r}_\diamond(\partial D)$ does not see an additive constant, we have
\begin{align}
\label{eq:trapezoidal}
\langle \tilde{f}_M - \tilde{f}, v \rangle_{\partial D} 
&= \sum_{m=0}^{M-1} \I_m^{\tilde{f}} v_0(\tilde{z}_m) - \int_{\partial D} \tilde{f} v_0 \, {\rm d} s \nonumber \\ 
&= \sum_{m=0}^{M-1}  \frac{2\pi}{M} \tilde{f}(\tilde{z}_m) v_0(\tilde{z}_m) - \int_{\partial D} \tilde{f} v_0 \, {\rm d} s. 
\end{align}
Observe that the right-hand side of \eqref{eq:trapezoidal} corresponds to the quadrature error of the trapezoidal rule for $\tilde{f} v_0$ over $\partial D$ on an equidistant grid. Hence, it follows from \cite[Theorem~1.1 \& Remark~1.2]{Kurganov09} that
$$
\big| \langle \tilde{f}_M - \tilde{f}, v \rangle_{\partial D} \big| \leq
\frac{C}{M^r} \left \| \frac{\partial^r (\tilde{f}v_0)}{\partial s^r} \right\|_{L^1(\partial D)},
$$
where $s$ denotes the arclength parameter on $\partial D$ and $C = C(r)>0$.
By expanding the derivative on the right-hand side and using the triangle and Schwarz inequalities, one easily sees that
$$
\big| \langle \tilde{f}_M - \tilde{f}, v \rangle_{\partial D} \big|
\leq  \frac{C'}{M^r} \| v_0 \|_{H^r(\partial D)} \| \tilde{f} \|_{H^r(\partial D)}, \qquad C' = C'(r) >0. 
$$
We now take the supremum over those $v \in H^r(\partial D)/\C$ that have unit norm in order to conclude that
\begin{equation}
\label{eq:bestesimate}
\big\| \tilde{f}_M - \tilde{f} \big\|_{H^{-r}(\partial D)} \leq  \frac{C'}{M^r} \| \tilde{f} \|_{H^r(\partial D)}
\end{equation}
for $r \in \N$.

Due to the boundedness of $\tilde{\RR}_{\tilde{\sigma}} - \tilde{\RR}_1: H_\diamond^r(\partial D) \to H^{q}(\partial D)/\C$ for any $q \in \R$, it follows from \eqref{eq:bestesimate} that
$$
\big\| \big(\tilde{\RR}_{\tilde{\sigma}} - \tilde{\RR}_1 \big)\big( \tilde{f}_M - \tilde{f}\big)\big\|_{H^q(\partial D)/\C} \leq \frac{C''}{M^r} \| \tilde{f} \|_{H^r(\partial D)} \qquad {\rm for} \ {\rm all} \ q \in \R
$$
and some $C''(q,r,\tilde{\sigma}) > 0$. The assertion now follows by combining \eqref{eq:Adef}, \eqref{eq:deltacompo} and the Sobolev embedding theorem.
\end{proof}

The optimal convergence rate in \eqref{eq:optconv} follows from the equidistant placement of the virtual point electrodes on $\partial D$, that is, from the efficiency of the trapezoidal rule for approximating integrals with periodic integrands~(cf.,~e.g.,~\cite{Kurganov09}). This was achieved by choosing the positions of the true electrodes on $\partial \Omega$ to be the conformal preimages of the the virtual equidistant ones. By comparing this observation,~e.g.,~with the M\"obius transform \eqref{mobius} used for magnifying a certain ROI in case $\Omega = D$, it becomes obvious that the true electrode pattern on $\partial \Omega$ is the densest close to the ROI.

\section{Simulating data for a half-plane}
\label{sec:halfplane}

Let $\Omega$ be the (unbounded) open upper half-plane 
$$
\{z\in\C \ | \ \im(z) > 0\} \, \simeq \, \{x\in \R^2 \ | \ x_2 > 0\}. 
$$
Although the general treatment of EIT in unbounded domains would require some extra considerations, here we exclusively restrict our attention to the following setting that is straightforward to handle and suits our purposes~\cite{Harhanen2010}:
\begin{itemize}
\item[I] The boundary current densities are (finite) linear combinations of delta distributions, i.e.,
\begin{equation}
\label{eq:deltacurr}
f_{\mathcal{I}} = \sum_{m=0}^{M-1} \mathcal{I}_m \delta_{z_m},
\end{equation}
where $\mathcal{I} \in \C_\diamond^M$ and $z_m$, $m=0,\dots, M-1$, are points on $\partial \Omega$,~i.e.,~on the horizontal axis.
\item[II] The conductivity $\sigma: \Omega \to \R_+$ is composed of a finite number of compactly supported $C^2$ inclusions with constant conductivity levels in unit background.
\end{itemize}
Under these assumptions, the Neumann problem \eqref{eq:neumann} has a unique solution $u \in H^1_{\rm loc}(\Omega)$ under the decay condition
$$
\lim_{|z|\to \infty} |u(z)| = 0.
$$
In particular, the relative boundary potentials corresponding to current patterns of the form \eqref{eq:deltacurr},~i.e.,
$$
(\RR_\sigma - \RR_1) f_{\mathcal{I}},\qquad \mathcal{I} \in \C_\diamond^M,
$$
are smooth and can be simulated using layer potentials. See \cite[Appendix]{Harhanen2010} for the details. In consequence, the relative current-to-potential map of the PEM remains well-defined for the half-space geometry.

What is more, the argumentation in the proof of Theorem~\ref{main_theorem} remains valid. That is, if (i) $\Psi$ is a conformal map of the unit disk $D$ onto the half-space $\Omega$,~i.e.,~a certain M\"obius transform, (ii) the points $\{ \tilde{z}_m \}_{m=0}^{M-1}$ form an uniform grid on $\partial D$ and (iii) $z_m = \Psi(\tilde{z}_m)$, $m = 0, \dots, M-1$, then the estimate \eqref{eq:optconv} is true for $\tilde{f} \in H^r_\diamond(\partial D)$ satisfying $[\tilde{f}]_M\in \C_\diamond^M$. (The only part in the proof of Theorem~\ref{main_theorem} that actually refers to $\Omega$ is \eqref{eq:deltacompo} which can be relatively easily validated for a half-space.)

When generating the data for the half-space geometry in our numerical examples, we first simulate the relative current-to-potential operator of the PEM \eqref{eq:Adef} by employing the boundary element code from~\cite{Harhanen2010}. Subsequently, we {\em approximate} the grid values of the (virtual) relative potential measurements
$$
(\tilde{\RR}_{\tilde{\sigma}} - \tilde{\RR}_1)  \tilde{f}
$$
for the Fourier current basis $\tilde{f} = \tilde{\phi}_n$, $1\leq |n| \leq N$, on the boundary of the virtual reconstruction domain $D$ with the help of \eqref{eq:optconv}. These can then be used for approximating the Fourier basis representation for the (relative) D-to-N map on the boundary of the virtual domain $\partial D$.

\bibliographystyle{plain}

\end{document}